\documentclass[12pt]{amsart}

\usepackage{verbatim}

\usepackage{amsthm}
\usepackage{amsfonts}
\usepackage{amsmath}
\usepackage{amssymb}
\usepackage{verbatim, color}

\usepackage{pictex}
\usepackage{url}

\usepackage{ragged2e} 

\begin{document}
 \newcounter{thlistctr}
 \newenvironment{thlist}{\
 \begin{list}%
 {\alph{thlistctr}}%
 {\setlength{\labelwidth}{2ex}%
 \setlength{\labelsep}{1ex}%
 \setlength{\leftmargin}{6ex}%
 \renewcommand{\makelabel}[1]{\makebox[\labelwidth][r]{\rm (##1)}}%
 \usecounter{thlistctr}}}%
 {\end{list}}

\thispagestyle{empty}

\newtheorem{Lemma}{\bf LEMMA}[section]
\newtheorem{Theorem}[Lemma]{\bf THEOREM}
\newtheorem{Claim}[Lemma]{\bf CLAIM}
\newtheorem{Corollary}[Lemma]{\bf COROLLARY}
\newtheorem{Proposition}[Lemma]{\bf PROPOSITION}
\newtheorem{Example}[Lemma]{\bf EXAMPLE}
\newtheorem{Fact}[Lemma]{\bf FACT}
\newtheorem{definition}[Lemma]{\bf DEFINITION}
\newtheorem{Notation}[Lemma]{\bf NOTATION}
\newtheorem{Remark}[Lemma]{\bf REMARK}
\newtheorem{theorem}[Lemma]{\bf THEOREM}

\newenvironment{Proof}{\noindent\bf Proof \rm}{$\hfill
\square$}

\newcommand{\restrict}{\mbox{$\mid\hspace{-1.1mm}\grave{}$}}
\newcommand{\covers}{\mbox{$>\hspace{-2.0mm}-{}$}}
\newcommand{\covered}{\mbox{$-\hspace{-2.0mm}<{}$}}
\newcommand{\notcover}{\mbox{$>\hspace{-2.0mm}\not -{}$}}

\newcommand{\boldalpha}{\mbox{\boldmath $\alpha$}}
\newcommand{\boldbeta}{\mbox{\boldmath $\beta$}}
\newcommand{\boldgamma}{\mbox{\boldmath $\gamma$}}
\newcommand{\boldxi}{\mbox{\boldmath $\xi$}}
\newcommand{\boldlambda}{\mbox{\boldmath $\lambda$}}
\newcommand{\boldmu}{\mbox{\boldmath $\mu$}}

\newcommand{\barzero}{\bar{0}}

\newcommand{\sfq}{{\sf q}}
\newcommand{\sfe}{{\sf e}}
\newcommand{\sfk}{{\sf k}}
\newcommand{\sfr}{{\sf r}}
\newcommand{\sfc}{{\sf c}}
\newcommand{\restr}{\negmedspace\upharpoonright\negmedspace}

\newcommand{\debaj}[2]{ #1 \to_H #2}

\newcounter{ecuacionDef} \setcounter{ecuacionDef}{0}
\renewcommand{\theecuacionDef}{\arabic{ecuacionDef}}

\newenvironment{ecuacionDef}[1]%
{
	\refstepcounter{ecuacionDef}
	\vspace{0.12cm}
	{\rm (E\theecuacionDef)}:  #1
	\vspace{0.12cm}
}%

\newcounter{numeroAxioma} \setcounter{numeroAxioma}{0}
\renewcommand{\thenumeroAxioma}{\arabic{numeroAxioma}}

\newenvironment{numeroAxioma}[1]%
{
	\refstepcounter{numeroAxioma}
	\vspace{0.12cm}
	(A\thenumeroAxioma):  #1
	\vspace{0.12cm}
}%

\title[Connexive Logics and 
semi-Heyting algebras]{Connexive logics  
and connexive semi-Heyting algebras}   
          
\dedicatory{Dedicated to Professor Kuznetsov}

\author{Juan M. CORNEJO* \and Hanamantagouda P. Sankappanavar}
\newcommand{\acr}{\newline\indent}

\address{\llap{*\,} Departamento de Matem\'atica\acr
Universidad Nacional del Sur\acr
Alem 1253, Bah\'ia Blanca, Argentina\acr
INMABB - CONICET}

\email{jmcornejo@uns.edu.ar}

\address{\llap{**\,}Department of Mathematics\acr
                              State University of New York\acr
                              New Paltz, New York, 12561\acr
                              U.S.A.}

\email{sankapph@hawkmail.newpaltz.edu}



\begin{abstract}
  
In this paper, we define and investigate a connexive logic, called ``Connexive semi-Heyting logic'' ($\mathcal{CSH}$ for short) and a new subvariety $\mathbb{CSH}$ of the variety $\mathbb{SH}$ of semi-Heyting algebras.  It is shown that the logic  $\mathcal{CSH}$ is implicative in the sense of Rasiowa, and is algebraizable with $\mathbb{CSH}$
as an equivalent algebraic semantics (in the sense of Blok and Pigozzi). 
 We also introduce the logics $\mathcal{AT}i$ and $\mathcal{BT}i$, $i=1,2$, along with the subvarieties $\mathbb{AT}i$ and $\mathbb{BT}i$, $i=1,2$ of $\mathbb{SH}$.  It is then shown that $\mathbb{AT}1 = \mathbb{AT}2 $  and $\mathbb{CSH} = \mathbb{BT}1 \subset \mathbb{BT}2 \subset \mathbb{AT}1$.
A 3-valued connexive semi-Heyting logic $\mathcal{CSH}3$ and its equivalent algebraic semantics $\mathbb{CSH}3$ are introduced and axiomatized; and it is  
then shown that $\mathcal{CSH}3$ is deductively equivalent to the 3-valued  intuitionistic logic.  
New characterizations of anti-Boolean semi-Heyting algebras are given.  We show that $\mathbb{BT}2 \cap \mathbb{SH}_c = \mathbb{V}(\mathbf{\bar{2}})$ and 
$\mathbb{SH}_c  \subset \mathbb{AT}1$, where $\mathbb{SH}_c $ is defined by
$x \to y \approx y \to x$.  It is  proved that the identity (AT1) is equivalent to the identity $x^* \to y^* \approx y^* \to x^*$ in $\mathbb{S}t\mathbb{SH}$ and also is equivalent to
$0\to 1 \approx 0$ in $\mathbb{SH}$.  We show that $\mathbb{AT}1 \cap \mathbb{EX} \subset \mathbb{BT}1$, where $\mathbb{EX}$ is defined by $x \to (y \to z) \approx y \to (x \to z)$.
The paper concludes with some further remarks, mentions some open problems for future research and proposes two new principles to be considered as Connexive Theses. 
\end{abstract}

\maketitle

\thispagestyle{empty}

\section{{\bf Introduction}} \label{SA}

Perhaps, it may not be an overstatement if we say that the heart and soul of a logic is its implication.
 There are several types of implications in Logic such as material implication, strict implication, relevance implication, etc., giving rise to different kinds of logics.  Connexive Logic deals with yet another implication, called connexive implication.  We noticed in April, 2020 that there are semi-Heyting algebras whose implication was connexive (see $\mathbf{\bar{2}}$ and $\mathbb{CSH}3$ in Section 2).
  Inspired by the theses of connexive logics, we introduce--and initiate a systematic investigation on-- the connexive semi-Heyting logics,   
connexive semi-Heyting algebras, and also display the interplay between them.  For some historical remarks  
on interplay between logics and (varieties of) algebras and for other such interplays, see \cite{Sa22b}.
 
 We develop, thanks to Rasiowa \cite{Ra74} and Blok and Pigozzi \cite{BlPi89}, an ``effective'' back-and-forth translation between connexive semi-Heyting logics and these new subvarieties of the variety $\mathbb{SH}$ of semi-Heyting algebras.  Using resources and tools from semi-Heyting algebras, we state and prove new theorems, which, in turn, will be translated into the theory of connexive logics.
 
 We need to start with some basic terminology of algebraic logic to fix the notation.  In this paper we consider the two languages $\langle \lor, \land, \to, \bot, \top \rangle$ and $\langle \lor, \land, \to, 0, 1 \rangle$ as essentially the same; however, we frequently use the former in the context of logics and the latter in the context of algebras. Let Var denote a fixed countable set of propositional variables. Let L  be a language.  The set of formulas (or terms) $Fm_L$ and the algebra $\mathbf{Fm_L}$ of formulas in the language L are inductively defined as usual. 
 
A {\it substitution-invariant consequence relation} $\vdash$ is a binary relation which is a subset of $P(Fm_L) \times Fm_L$ satisfying the following conditions, for all $\Gamma$, $\Delta \subseteq Fm_L$, $\phi$, $\psi \in Fm_L$, and $\sigma$ an L-substitution (an endomorphism of the formula algebra $\mathbf{Fm}_L$):

(1) $\Gamma \vdash \phi$ whenever $\phi \in \Gamma$ \quad  (Reflexivity);

(2) If $\Gamma \vdash \phi$ and $\Gamma \subseteq \Delta$, then $\Delta \vdash  \phi$ \quad (Monotonicity);

(3) If $\Delta \vdash  \phi$ and $\Gamma \vdash \psi$ for every $\psi \in \Delta$, then $\Gamma \vdash \phi$ \quad (Cut);

(4) If $\Gamma \vdash \phi$, then $\sigma(\Gamma) \vdash \sigma(\phi)$. \quad (Substitution-invariance).
 
 A \emph{(propositional) logic} $\mathcal L$ is a pair $\mathcal L =   \langle L, \vdash_\mathcal{L} \rangle$, where $L$ is a language and 
 $\vdash_\mathcal{L}$ is a substitution-invariant consequence relation.   
We will present logics by means of their ``Hilbert style'' axioms and inference rules.
 
\subsection{Brief Historical Remarks on Connexive Logic}

\

Connexive Logic has its roots in ancient Greek philosophy.
It is frequently mentioned that Aristotle believed in the validity of statements like ``A proposition cannot be implied by its own negation.''  This statement has come to be known as \emph{Aristotle's thesis}, which is clearly false in the classical logic. 

Boethius, a medieval philosopher, believed in the validity of statements like ``no proposition entails both a proposition and its negation."  This statement has come to be known as ``Boethius’ thesis.``

In 1963, Storrs McCall, in his Ph.D. Thesis \cite{Mc63}, writes: 
 
\begin{quote}
``But there exists also a type of implication, intuitively plausible, which is nonclassical not only in being more restrictive, but in satisfying certain
theses which are classically false. These theses are exceedingly
venerable, dating back to Aristotle and Boethius, but, despite their
plausibility, have been generally rejected by logicians since. It has
not been noticed, however, that in Sextus Empiricus reference is made
to a species of Stoic implication which fits them perfectly.''
\end{quote}

He continues to say:

\begin{quote}
``In this work formal recognition is given to this species of
implication, known as connexive implication. It is shown that none of
the well-known systems of prepositional logic is connexive.'' 
\end{quote}
The name ``Connexive Logic'' originates from McCall's Ph.D. thesis. 
 McCall adopted the term ``connexive'' from Bochenski's book \cite{Bo61}.
Storrs McCall is generally regarded as the founder of Modern Connexive Logic.  We should also mention here that
the recent analysis of E.J.
Nelson’s connexive logic by Mares and Paoli  
seems to suggest that E.J. Nelson’s work may have been a precursor to McCall's work.

 Omori and Wansing, the editors of the recently published book, entitled ``60 Years of Connexive Logic'' write:

\begin{quote}
``Unlike the most well known systems of non-classical logic, systems of connexive logic are contra-classical in that they both reject certain classically valid principles and validate schemata that are not valid classically. The history of modern formal connexive logic may be seen to have started with Storrs McCall’s dissertation “Non-classical Propositional Calculi" (Oxford) 1963, thus roughly 60 years ago. While at the turn of the 21st century connexive logic was a virtually dead research program, the situation has changed significantly after the inclusion of an entry on connexive logic in the Stanford Encyclopedia of Philosophy in 2006 and the beginning of a series of annual workshops on connexive logic in 2015. Nowadays, connexive logic is a vibrant area, and the present volume offers an exciting glimpse on recent work in connexive logic. $\cdots$.''  
\end{quote}

For more details on the motivation, the origin and the history of connexive logics, see \cite{Mc12} and \cite{Wa22}. 

We will conclude our brief historical remarks on Connexive Logic by mentioning that we are not aware--but are very much interested in knowing-- of any unbiased, modern investigation into the question of whether Indian philosophy or Chinese philosophy had also considered statements akin to the Aristotle's  and Boethius' Theses.

\subsection{Connexive logics: Formally Speaking}
\

It is generally agreed that Connexive Logic is defined to be a logic that has Aristotle's and Boethius' Theses as theorems.  We will now present these theses more formally.  In the formal presentation, each of these theses gives rise to two theses.
\begin{definition}
 A logic $\mathcal L$ in a language {\bf L}, (possibly) extending the language $\langle \to, \neg \rangle$  is a \emph{connexive logic}  
 if the following Aristotle's Theses and Boethius' Theses are theorems in $\mathcal L$:\\

\noindent {\bf Aristotle's Theses}: 

\indent  \rm(AT1\rm)  \quad $\neg(\neg \alpha \to \alpha)$, 

\indent  (AT2)  \quad $\neg(\alpha \to \neg \alpha)$,\\

\noindent {\bf Boethius' Theses}:

\indent (BT1)  \quad $(\alpha \to \beta) \to \neg(\alpha \to \neg \beta)$, 
  
\indent  (BT2) \quad $(\alpha \to \neg \beta) \to \neg(\alpha \to \beta).$\\

A connexive logic $\mathcal L$ is \emph{non-symmetric} if the implication is non-symmetric, i.e., $(\phi \to \psi) \to(\psi \to \phi)$
 fails to be a theorem (so that $\to$ can not be understood as a bi-conditional).
\end{definition}

Some logicians, like 
Kapsner (see \cite{Ka12}) consider even more stringent requirements and say that a connexive logic is \emph{strongly connexive} if it satisfies the two additional conditions:

(a) In no model, $\alpha \to \neg\alpha$ is satisfiable ($\alpha$ being an arbitrary formula);

(b) In no model, $\alpha \to \beta$ and $\alpha \to \neg\beta$ are simultaneously satisfiable ($\alpha$ and $\beta$ are arbitrary formulas).

In this paper{\footnote {Results of this paper were presented by the second author in an invited lecture at AAA107--Workshop on General Algebra, Bern, Switzerland,
 June 20-22, 2025.}, we define and investigate a non-symmetric, (strongly) connexive logic, called ``Connexive Semi-Heyting Logic ($\mathcal{CSH}$)'' that is implicative in the sense of Rasiowa and hence, it is algebraizable in the sense of Blok and Pigozzi. 

Jarmuzek and Malinowski \cite{JaMa19} have introduced the notion of a ``quasi-connexive'' logic.  

A logic is \emph{quasi-connexive} if it is not connexive, but at least one of (AT1), (AT2), (BT1) and (BT2) is a theorem in the logic.

Motivated by Aristotle's and Boethius' theses, we go a step further and we introduce logics $\mathcal{AT}1$, $\mathcal{AT}2$, $\mathcal{BT}1$ and $\mathcal{BT}2$ and examine their mutual relationships
 and their relation with $\mathcal{CSH}$.

\begin{definition}
Let $\mathcal{L}$ be a logic.  Then $\mathcal{L}$ is  \rm(AT1)-\emph{connexive} if the formula (AT1) is provable in $\mathcal{L}$.  Similarly, (AT2)-\emph{connexive}, (BT1)-\emph{connexive}, 
and (BT2)-\emph{connexive} logics are defined.  $\mathcal{L}$ is {\em anti-Boolean} (or {\em contra-classical}) if the formulas $\neg \neg \alpha \leftrightarrow \alpha$ and $(\bot \to \top) \leftrightarrow \bot$ are provable in $\mathcal{L}$.  $\mathcal{L}$ is {\em anti-intuitionistic} (or {\em anti-Heyting}) if the formula  $(\bot \to \top) \leftrightarrow \bot$ is provable in $\mathcal{L}$, where $\alpha \leftrightarrow \beta := (\alpha \to \beta) \land  (\beta \to \alpha)$.                                                               
\end{definition}

\subsection{Semi-Heyting Algebras}

\

Semi-Heyting algebras were introduced by the second author in 1985 as an abstraction of Heyting algebras and were published much later in 2007 (\cite{Sa07}).  See also 
\cite{Ab11, Ab11a, AbCoVa13, CoMoSaVi20, CoSa19}).

An algebra $\mathbf{A} = \langle L; \land, \lor, \to, 0, 1\rangle$ is a semi-Heyting algebra if the following conditions hold:
\begin{enumerate}
	\item[(SH1)] $\langle L; \land, \lor, 0,1\rangle$ is a lattice with $0$ and $1$,
	\item[(SH2)] $x \land (x \to y) \approx x \land y$,
	\item[(SH3)] $x \land (y \to z) \approx x \land [(x \land y) \to (x \land z)]$,
	\item[(SH4)] $x \to x \approx 1$.
\end{enumerate}
A semi-Heyting algebra is a Heyting algebra if it satisfies:
$$(x \land y) \to x \approx 1.$$
We will denote the variety of semi-Heyting algebras by $\mathbb{SH}$ and that of Heyting algebras by $\mathbb{H}.$

Semi-Heyting algebras share (see \cite{Sa07}) some rather strong properties with Heyting algebras.   
 For example, (1) every interval in a semi-Heyting algebra is also pseudocomplemented, (2) the variety SH is arithmetical, (3) the variety SH has EDPC (equationally definable principal congruences), and (4) Congruences on the algebras in V are determined by filters.   Moreover, there is a rich supply of algebras in SH. For example, It is known that there are in SH, up to isomorphism, two 2-element algebras, ten 3-element algebras, only one of which, of course, is a Heyting algebra and 160 algebras on a 4-element chain (see \cite{Sa07, Ab11,  Ab11a, AbCoVa13, CoMoSaVi20, CoSa19}).

\subsection{Connexive Semi-Heyting Logic $\mathcal{CSH}$}  
\

\medskip 

In \cite{Co11}, the first author has introduced “semi-intuitionistic logic” (also called, semi-Heyting logic''), which is algebrizable (in fact, implicative, in the sense of Rasiowa) with the variety of semi-Heyting algebras as its equivalent algebraic semantics (see  Section 3).

    The fact that the identity $0 \to 1 \approx 0$ holds in (infintely) many semi-Heyting algebras led us to consider, in 2020, the possibility that there might be connexive logics arising from semi-Heyting algebras.  We noticed in May 2020 that it was indeed the case. The matrices associated with algebras $\bar{\mathbf{2}}$ and $\mathbf{CSH3}$ (named $\mathbf{L_9}$ in \cite{CoSa22a}) were our first such examples, while (AT1) and (AT2) were true in $\mathbf{L_{10}}$.   In fact, we wrote the following paragraph in \cite{CoSa22a}:

``Many of the extensions of the logics $\mathcal{SH}$ and $\mathcal{DHMSH}$, to our surprise, turn out to be connexive logics with $\neg \alpha := \alpha \to \bot$.  We present a few of these below.  (More will be said in a paper which is in preparation.)''  

It is the present paper and \cite{CoSa25} that were referred to in the parenthetical statement of the  preceding paragraph. 
In this paper we will present, with proofs, the various relationships that exist among Aristotle's Theses and Boethius' Theses in the logic $\mathcal{SH}$ that
were promised in the paper \cite{CoSa22a}.  

The paper is organized as follows: Section 2 contains definitions, notation and some new arithmetical properties of semi-Heyting algebras that are needed later in the paper (and also possibly for future research).
Section 3 introduces connexive semi-Heyting logic $\mathcal{CSH}$ as an extension of the logic $\mathcal{SH}$, 
which is one of the two central notions of this paper.  It also introduces new logics $\mathcal{AT}1$, $\mathcal{AT}2$, $\mathcal{BT}1$, and $\mathcal{BT}2$, as axiomatic extensions of $\mathcal{SH}$.  As such, these logics are all implicative in the sense of Rasiowa. 
It is then shown that $\mathcal{CSH}$ is implicative and hence algebraizable.  In Section 4, we introduce a new subvariety $\mathbb{CSH}$ of $\mathbb{SH}$ consisting of ``connexive semi-Heyting algebras'', which is the other central notion of this paper.  The varieties $\mathbb{AT}1$, $\mathbb{AT}2$, $\mathbb{BT}1$, and $\mathbb{BT}2$ are also introduced.  It is shown that the varieties $\mathbb{CSH}$, $\mathbb{AT}1$, $\mathbb{AT}2$, $\mathbb{BT}1$, and $\mathbb{BT}2$ are respectively the equivalent algebraic semantics of the logics $\mathcal{CSH}$, $\mathcal{AT}1$, $\mathcal{AT}2$, $\mathcal{BT}1$, and $\mathcal{BT}2$.
Section 5 exhibits mutual relationships among AT1, AT2, BT1, and BT2, and $\mathcal{CSH}$ relative to the semi-Heyting logic $\mathcal{SH}$, by showing that $\mathbb{BT}1 = \mathbb{CSH}$, $\mathbb{AT}1= \mathbb{AT}2$ and $\mathbb{BT}1 \subset \mathbb{BT}2 \subset \mathbb{AT}1$.
In Section 6 we axiomatize the variety $\mathbf{CSH3}$ generated by the 3-valued connexive semi-Heyting algebra, and by translation, we obtain an axiomatization of the logic $\mathcal{CSH}3$.  It is then shown that the logic $\mathcal{CSH}3$ is deductively equivalent to the the 3-valued Heyting logic $\mathcal{H}3$.
In Section 7, we give new characterizations of anti-Boolean semi-Heyting algebras.   Furthermore, we show that $\mathbb{BT}2 \cap \mathbb{SH}_c = \mathbb{V}(\mathbf{\bar{2}})$ and 
$\mathbb{SH}_c  \subset \mathbb{AT}1$, where $\mathbb{SH}_c $ denotes the variety of commutative (i.e., satisfying the identity: $x \to y \approx y \to x$)  semi-Heyting algebras.  Section 8 proves that the identity (AT1) and the weak commutative identity $x^* \to y^* \approx y^* \to x^*$ are equivalent in the variety $\mathbb{S}t\mathbb{SH}$ of Stone semi-Heyting algebras, while in Section 9, it is proved that the identity (AT1) and the anti-Heyting identity $0\to 1 \approx 0$ are equivalent in $\mathbb{SH}$.  
In Section 10, we show that $\mathbb{AT}1 \cap \mathbb{EX} \subset \mathbb{BT}1$, where  
$\mathbb{EX}$ is the subvariety of $\mathbb{SH}$ defined by the exchange identity $x \to (y \to z) \approx y \to (x \to z)$. 
Section 11 makes some further remarks, mentions some open problems for further research and proposes two new principles as Connexive Theses.

\medskip

\vspace{1cm}
\section{\bf {Preliminaries}} \label{SB}

In this section we recall some notions and results needed to make this paper as self-contained as possible.  However, for more basic information, we refer the reader to \cite{BuSa81} for concepts and results in universal algebra, to \cite{BaDw74} for distributive lattices, and for algebraic logic to \cite{Ra74} and \cite{Fo17}.  We then develop some new arithmetical properties about semi-Heyting algebras which will be useful later in this paper and in future research.

We should mention that, in the rest of the paper, we use $x, y, z$ both formally to represent  variables in identities, and mathematically to denote elements of an algebra.  This, we believe, would not cause confusion, since the context will make their meaning clear.  However, we use $a,b,c$ only to denote elements of an algebra.

\subsection{Some important examples of $\mathcal{SH}$} 

\  

 The following algebras, 
 $\mathbf{2}$ and  $\mathbf{\bar{2}}$
 are (with $0 <1$), up to isomorphism, the only two 2-element algebras in $\mathbb{SH}$. \\

$\mathbf{2}$:
\begin{tabular}{r|rr}
	$\to$: & 0 & 1\\
	\hline
	0 & 1 & 1 \\
	1 & 0 & 1
\end{tabular} \hspace{5cm}
$\mathbf{\bar{2}}$: 
\begin{tabular}{r|rr}
	$\to$: & 0 & 1\\
	\hline
	0 & 1 & 0 \\
	1 & 0 & 1
\end{tabular} \hspace{.5cm}
\smallskip
\begin{center} Figure 1
\end{center}

\indent It was shown in \cite{Sa07} 
that there are, up to isomorphism, ten $3$-element 
semi-Heyting algebras, We will list some of them which will be useful in this paper, below-- with their $\to$ operations only-- in figure $2$ below, where $0 < a < 1$:

\vspace{1cm}

\setlength{\unitlength}{1cm} \vspace*{.5cm}
\begin{picture}(15,2)
\put(1.5,-.4){\circle*{.15}}          

\put(1.5,1){\circle*{.15}}

\put(1.5,2.3){\circle*{.15}}

\put(1.1,-.4){$0$}                      

\put(1.1,.9){$a$}

\put(1.1,2.1){$1$}

\put(.1,.7){${\bf L}_1$ \ :}

\put(1.5,-.4){\line(0,1){2.7}}

\put(2.2,1){\begin{tabular}{c|ccc}
 $\to$ \ & \ 0 \ & \ $a$ \ & \ 1 \\ \hline
  & & & \\
  0 & 1 & 1 & 1 \\
  & & & \\
  $a$ & 0 & 1 & 1 \\
  & & & \\
  1 & 0 & $a$ & 1
\end{tabular}
}
\end{picture}

\medskip
\setlength{\unitlength}{1cm} \vspace*{1.5cm}
\begin{picture}(15,2)
\put(1.5,-.4){\circle*{.15}}

\put(1.5,1){\circle*{.15}}

\put(1.5,2.3){\circle*{.15}}

\put(1.1,-.4){$0$}

\put(1.1,.9){$a$}

\put(1.1,2.1){$1$}

\put(-.2, 1.3){${\bf CSH3:}$\:} 

\put(1.5,-.4){\line(0,1){2.7}}

\put(2.2,1){\begin{tabular}{c|ccc}
 $\to$ \ & \ 0 \ & \ $a$ \ & \ 1 \\ \hline
  & & & \\
  0 & 1 & {\color{red}0} & {\color{red}0} \\
  & & & \\
  $a$ & 0 & 1 & 1 \\
  & & & \\
  1 & 0 & $a$ & 1
\end{tabular}
}

\put(8.2,-.4){\circle*{.15}}

\put(8.2,1){\circle*{.15}}

\put(8.2,2.3){\circle*{.15}}

\put(7.8,-.4){$0$}

\put(7.8,.9){$a$}

\put(7.8,2.1){$1$}

\put(6.7,.7){${\bf L}_{10}$ \ :}

\put(8.2,-.4){\line(0,1){2.7}}

\put(8.9,1){\begin{tabular}{c|ccc}
 $\to$ \ & \ 0 \ & \ $a$ \ & \ 1 \\ \hline
  & & & \\
  0 & 1 & {\color{red}0} & {\color{red}0} \\
  & & & \\
  $a$ & 0 & 1 & {\color{red}$a$} \\
  & & & \\
  1 & 0 & $a$ & 1
\end{tabular}
}

\put(6,-1.4){Figure 2} \label{page_algebra2}  
\end{picture}

\vspace{1.5cm}
(Note that $\mathbf{L_1}$ and $\mathbf{L_{10}}$ are names borrowed from \cite{Sa07}, 
while  $\mathbf{CSH3}$ was named $\mathbf{L_9}$ in \cite{Sa07}.)

\vspace{1.5cm}

\subsection{Arithmetical Properties of Semi-Heyting Algebras}

\

\medskip
Here we present a few arithmetical properties of semi-Heyting algebrras that will be needed in the sequel and for future work.

The proof of the following lemma is straightforward.

\begin{Lemma} \label{2.2}
Let ${\mathbf L} \in \mathbb{SH}$ and let $x,y, z \in L$.  Then
\begin{enumerate}
 \item[{\rm(i)}] $1 \to x =x$,
\item[{\rm(ii)}] $x  \leq (x  \lor y) \to x$, \label{B.63}
\item[{\rm(iii)}] $x  \leq x \to 1$, \label{C.170}
\item[{\rm(iv)}] $x \land [(x \to y) \to z] = x \land (y \to  z)$,  
\item[{\rm(v)}] $x \land [y \lor (x \to z)] = x \land (y \lor z)$, 
\item[{\rm(vi)}]  $x \land [(x \lor y) \to z] = x \land z$,  

\item[{\rm(vii)}]  $x \land [(x \land y) \to z] = x \land (y \to z)$,  

\item[(viii)] \label{Lemma_3.70} 

$x \land [(x \to y) \lor z] = x \land (y \lor z)$,   

\item[(ix)] \label{Lemma_4.62} 
$y \land [x \to (y \land z)] = y \land (x \to  z)$, 
 
\item[(x)] \label{Lemma_4.75}  
$x \to (y  \land z) \geq y \land (x \to  z)$, 

\item[(xi)] \label{Lemma_4.93} 
$x\  \leq \ y \to (x \land y)$.  
\end{enumerate}
\end{Lemma}

The following lemma holds in semi-Heyting algebras as they are pseudocomplemented.

\begin{Lemma}
Let ${\mathbf L} \in \mathbb{SH}$. Then ${\mathbf L}$ satisfies the following identities:
\begin{enumerate}
\item  $(x^{**} \land y)^*= (x^{**} \land y^{**})^* = (x \land y)^*$. \label{1245}
\item  $x \land (x^* \land y)^* = x.$ \label{147}
\end{enumerate}
\end{Lemma}

The following lemma is easy to verify using (SH3) from the definition and elementary properties
of semi-Heyting algebras. 

\begin{Lemma} \label{634} \label{659.1}  \label{60}  \label{866} \label{E37}  \label{232}
	Let ${\mathbf L} \in \mathbb{SH}$. Then ${\mathbf L}$ satisfies the following identities:
	\begin{enumerate}
		\item $x \land [\{(x \land y) \to (x \land z)\} \to u] = x \land [(y \to z) \to u]. $ 
		\item   $x \land [\{y \to (x \land z)\} \to u] = x \land ((y \to z) \to u). $ 
		\item  $[(x \to y) \to z] \leq [\{(x \land (y \to z)) \to y\} \to z].$ 
		\item  $x \land (y \to (x \to z)) = x \land (y \to z).$ 
		\item  	$x \land [y \to ((x \land z) \to (x \land u))] = x \land [y \to (z \to u)].$ 
		\item   $x \land ((x \to y) \to z) = x \land (y \to z).$  
	\end{enumerate}
\end{Lemma}

\vspace{1.1cm}

\section{Connexive Semi-Heyting logic $\mathcal{CSH}$ and the logics $\mathcal{AT}1$, $\mathcal{AT}2$, $\mathcal{BT}1$,  and $\mathcal{BT}2$} 

\

In this section we wish to 
define and initiate the investigation into the connexive logic, and related logics $\mathcal{AT}i$ and $\mathcal{BT}i$, $i=1,2$.  
Since they are going to be certain axiomatic extensions of the semi-Heyting logic $\mathcal{SH}$ (also called, semi-ntuitionistic logic) and inherit the properties of implicativeness and algebraizability of the logic $\mathcal{SH}$, it is imperative that   
we first recall some definitions and results of the Abstract Algebraic Logic \cite{BlPi89, Fo17} and also from \cite{Sa07, CoSa22a}.

\subsection{The Semi-Heyting logic $\mathcal{SH}$}

\

\begin{definition} {\rm  \cite{CoVi15}}
The semi-intuitionistic logic $\mathcal{SH}$ (also called $\mathcal{SI}$) is defined in the language 
$\{\vee, \land,  \to,  \bot, \top \}$ and it has the following   
axioms and the inference rule:\\

{\bf  \noindent  AXIOMS:}

\noindent \rm \numeroAxioma{$\debaj{\alpha}{(\alpha \vee \beta)}$,} \label{axioma_supremo_izq}

\noindent \numeroAxioma{$\debaj{\beta}{(\alpha \vee \beta)}$,} \label{axioma_supremo_der}

\noindent \numeroAxioma{$\debaj{(\debaj{\alpha}{\gamma})} {  [\debaj{(\debaj{\beta}{\gamma})} {(\debaj{(\alpha \vee \beta)}{\gamma})}]}$,} \label{axioma_supremo_cota_inferior}

\noindent \numeroAxioma{$\debaj{(\alpha \wedge \beta)}{\alpha}$,} \label{axioma_infimo_izq}

\noindent \numeroAxioma{$\debaj{(\debaj{\gamma}{\alpha})} { [\debaj{(\debaj{\gamma}{\beta})}{(\debaj{\gamma}{(\alpha \wedge \beta)})}]}$,} \label{axioma_infimo_cota_superior}

\noindent \numeroAxioma{$\top$,} \label{axioma_top}

\noindent \numeroAxioma{$\debaj{\bot}{\alpha}$,} \label{axioma_bot}

\noindent \numeroAxioma{$\debaj{(\debaj{(\alpha \wedge \beta)}{\gamma})}{(\debaj{\alpha}{(\debaj{\beta}{\gamma})})}$,} \label{axioma_condicRes_InfAImplic}

\noindent \numeroAxioma{$\debaj{(\debaj{\alpha}{(\debaj{\beta}{\gamma})})}{(\debaj{(\alpha \wedge \beta)}{\gamma})}$,} \label{axioma_condicRes_ImplicAInf}

\noindent \numeroAxioma{$\debaj{(\debaj{\alpha}{\beta})}{(\debaj{(\debaj{\beta}{\alpha})}{(\debaj{(\alpha \to \gamma)}{(\beta \to \gamma)})})}$,} \label{axioma_BuenaDefImplic2} \label{axioma_ImplicADerecha}

\noindent \numeroAxioma{$\debaj{(\debaj{\alpha}{\beta})}{(\debaj{(\debaj{\beta}{\alpha})}{(\debaj{(\gamma \to \beta)}{(\gamma \to \alpha)})})}$.} \label{axioma_BuenaDefImplic1} \label{axioma_ImplicAIzquierda} \\

  {\bf RULE OF INFERENCE:}
 \smallskip
	
\qquad (SMP):   From $\phi$ and $\phi \to_H {\gamma}$, deduce $\gamma$  ({\rm semi-Modus Ponens}).\\

\end{definition}

\subsection{Implicativeness of the logic $\mathcal{SH}$ } 

\

\smallskip

An important class of logics called “standard systems of implicative extensional propositional calculus” and a class of algebras, associated with each of those classes of logics were introduced by Rasiowa \cite[page 179]{Ra74}, through a generalization of the classical Lindenbaum-Tarski process.  We will refer to these logics as “implicative logics in the sense of Rasiowa” (“implicative logics”, for short). These logics have played a pivotal role in the development of Abstract Algebraic Logic. We now recall the definition of implicative logics. Our presentation follows Font \cite{Fo17}.

\begin{definition} \cite{Ra74} 
	Let $\mathcal L$  be a logic in a language that includes a binary connective $\to$, either primitive or defined by a term in exactly two variables. Then $\mathcal L$ is called an implicative logic with respect to the binary connective $\to$ if the following conditions are satisfied:
	\begin{itemize}
		\item[\rm (IL1)] $\vdash_{\mathcal L} \alpha \to \alpha$. 
		\item[\rm (IL2)] $\alpha \to \beta, \ \beta\to \gamma \vdash_{\mathcal L} \alpha \to \gamma$.
		\item[\rm (IL3)] For each connective $f$ of arity $n>0$ in the language of $\mathcal{L}$,\\ 
		$\left\{\begin{array}{c}
		\alpha_1 \to \beta_1, \ldots, \alpha_n \to \beta_n, \\
		\beta_1 \to \alpha_1, \ldots, \beta_n \to \alpha_n
		\end{array}
		\right\} \vdash_{\mathcal L} f(\alpha_1, \ldots, \alpha_n) \to f(\beta_1, \ldots, \beta_n)$.
		\item[\rm (IL4)] $\alpha, \alpha \to \beta \vdash_{\mathcal L} \beta$.
		\item[\rm (IL5)] $\alpha \vdash_{\mathcal L} \beta \to \alpha$.
	\end{itemize}
\end{definition}

\begin{definition} {\rm \cite[Definition 6, page 181]{Ra74}} \label{definicion_lalg}
	Let $\mathcal L$ be an implicative logic in the language $L$  with an implication connective $\to$, either primitive or defined by a term in exactly two variables.  An $\mathcal L$-algebra is an algebra $\mathbf A$ in the language $L$ that has an element $\top$ with the following properties:
	\begin{quote}
		\begin{itemize}
			\item[{\rm (LALG1)}] For all $\Gamma \cup \{\phi\} \subseteq \mbox{\it Fm} $ and all $h \in Hom({\bf Fm}_{\bf L}, \mathbf A)$, if $\Gamma \vdash_{\mathcal L} \phi$ and $h \Gamma \subseteq \{\top\}$ then $h \phi = \top$, 
			\item[{\rm (LALG2)}] For all $a,b \in A$, if $a \to b = \top$ and $b \to a = \top$ then $a=b$.
		\end{itemize}
	\end{quote}
	The class of $\mathcal L$-algebras is denoted by $\mathbf{Alg^*\mathcal{L}}$.
\end{definition}

The following Theorem was proved in \cite{CoVi15} by showing that $\mathbf{Alg^*\mathcal{SH}}= \mathbb{SH}$.

\begin{theorem} \cite{CoVi15} \label{teorema_SH_implicativa} 
	The logic $\mathcal{SH}$ is implicative with respect to the connective $\to_H$.
\end{theorem}

\subsection{Algebraizability of the logic $\mathcal{SH}$ } 

\

Here we first recall some relevant notions and results from Abstract Algebraic Logic  
(see \cite[Section 2.1]{BlPi89}, \cite{FJP03}, or \cite{Fo17}).

Let $L$ denote a language. 
 Identities in $L$ are ordered pairs of $L$-formulas that will be written in the form $\alpha \approx \beta$.  The set of identities in L is denoted by $Eq_L$.  
An interpretation $h$ in $\mathbf A$ is a homomorphism from the algebra of formulas $\mathbf{Fm}_L$ into $\mathbf{A}$.  An interpretation $h$ in $\mathbf A$
satisfies an identity $\alpha \approx \beta$ if $h \alpha = h \beta$.  
We denote this satisfaction relation by the notation: $\mathbf{A} \models_h 
\alpha \approx \beta $. An algebra $\mathbf A$ {\it satisfies the equation} $\alpha \approx \beta$ if all the interpretations in $\mathbf A$ satisfy it; in symbols,
$$\mathbf{A} \models \alpha \approx \beta \mbox{ if and only if } \mathbf{A} \models_h  
\alpha \approx \beta, \mbox{ for all interpretations in } \mathbf A.$$
A class $\mathbb K$ of algebras  {\it satisfies the identity} $\alpha \approx \beta$ when all algebras in $\mathbb K$ satisfy it; i.e.
$$\mathbb{K} \models \alpha \approx \beta \mbox{ if and only if } \mathbf{A}\models \alpha \approx \beta, \mbox{ for all }\mathbf A \in \mathbb K.$$ 
  
If $\bar{x}$ is a sequence of variables and $h$ is an interpretation in $\mathbf{A}$, then we write $\bar{a}$ for $h(\bar{x})$.
For a class $\mathbb K$ of $\mathbf L$-algebras, 
we define the relation $\models_{\mathbb K}$ that holds between a set $\Delta$ of identities and a single identity $\alpha \approx  \beta$ as follows: \\ 
$$\Delta \models_{\mathbb K} \alpha \approx \beta \text{ if and only if }$$

for every $\mathbf A \in \mathbb{K}$ and every interpretation $\bar{a}$ of the variables of $\Delta \cup  \{\alpha \approx \beta\}$ in $\mathbf{A}$, 

if $\phi^{\mathbf A}(\bar{a} ) = \psi^{\mathbf A} (\bar{a}) $,  for every $\phi \approx \psi \in \Delta$, then \  
$\alpha^{\mathbf A}(\bar{a}) = \beta^{\mathbf A}(\bar{a})$.\\

In this case, we say that $\alpha \approx \beta$ is a $\mathbb K$-consequence of $\Delta$. 
The relation $\models_K$ is called the {\it semantic equational consequence relation} determined by $\mathbb{K}$. \\

\begin{definition} {\rm \cite[Definition 2.2]{BlPi89} {\rm \cite[Definition 3.4]{Fo17}}}
Given a language $L$ and a logic $\mathcal{L} =\langle L, \vdash_{\mathcal L}\rangle$, a class $\mathbb K$ of  L-algebras is called an ``algebraic semantics'' for $\mathcal{L}$,   if \ $\vdash_{\mathcal L}$ can be interpreted in $\models_\mathbb K$ in the following sense:

 There exists a finite set $\{\delta_i (x) \approx \epsilon_i (x) : i < n\}$ 
 of identities with a single variable $x$ such that, for all $\Gamma \cup {\phi} \subseteq  Fm$ and each $j < n$,
$$
 \qquad \Gamma \vdash_{\mathcal L} \phi \ \Leftrightarrow \ \{\delta_i(\psi) \approx \epsilon_i(\psi) : i < n,  \psi  \in \Gamma\} \models_{\mathbb{K}}  \delta_j(\phi) \approx \epsilon_j(\phi).$$
The identities $\delta_i \approx  \epsilon_i$, for $i < n$, satisfying the above conidtion are called ``defining identities'' for 
$\mathcal L$ and $\mathbb K$.\\
\end{definition}

\begin{definition} {\rm \cite[Definition 2.8]{BlPi89}}, {\rm \cite[Definition 3.11]{Fo17}} \label{def2.8}
Let $\mathbb K$ be an algebraic semantics for a logic
$\mathcal L$  with defining identities $\delta_i \approx  \epsilon_i$, for $i < n$. 
Then $\mathbb K$ is said to be ``equivalent to'' $\mathcal L$ if there exists a finite set $\Delta_j(p,q)$, for $j < m$, of formulas with two variables $p, q$ such that for every identity $\phi  \approx \psi$, for $i<n$, and for $j<m$, 
 $$ \qquad   \phi  \approx \psi \models_{\mathbb K}  \{\delta_i(\Delta_j(\phi, \psi)) \approx \epsilon_i(\Delta_j(\phi, \psi)) : j<m, i<n,\}$$
 \begin{center} \text{ and  }
 \end{center}
 $$ \{\delta_i(\Delta_j(\phi, \psi)) \approx \epsilon_i(\Delta_j(\phi, \psi)) : j<m, i<n,\} \models_{\mathbb{K}}  \phi  \approx \psi.$$
 
 The set $\{\Delta_j(p,q) : j < m\}$ of formulas with two variables, satisfying the above condition is called a set of ``equivalence formulas'' for $\mathcal L$ and $\mathbb K$.
 A logic $\mathcal L$ is said to be ``algebraizable'' if it has an equivalent algebraic semantics $\mathbb K$.
\end{definition}

The following theorem, proved in \cite{BlPi89}, is crucial in what follows. 

\begin{Theorem} {\rm(\cite{BlPi89}, \cite[Proposition 3.15]{Fo17}}) \label{EquivalentSemanticsTheorem}
Every implicative logic  $\mathcal L$ is algebraizable with respect to the class $Alg^*\mathcal L$ and the algebraizability is witnessed by the defining identity $x \approx x \to x$ (or equivalently, $x \approx 1$) and the equivalence formulas 
$\Delta =\{p \to q, q \to p\}$. 
\end{Theorem}

The following result is taken from \cite[Corollary 5.5]{CoSa22a} which is crucial in the sequel. 
 
\begin{Corollary} \label{CorExt} 
 The logic $\mathcal{SH}$ is algebraizable, and  the variety  
 $\mathbb{SH}$ is the equivalent algebraic semantics for $\mathcal{SH}$ with the defining identity $x \approx x \to_H x$ \rm(or equivalently, $x \approx 1$) and the equivalence formulas $\Delta = \{p \to_H q, q \to_H p \}$.
 \end{Corollary}

\subsection{Axiomatic Extensions  of the logic $\mathcal{SH}$} 

\

An  \emph{axiomatic extension} of a logic $\mathcal{L}$ is a logic obtained by adjoining new axioms to those of $\mathcal{SH}$ but keeping the rules of inference the same as in $\mathcal L$.  

 In what follows, we frequently use the term ``extension'' for ``axiomatic extension''. 
 It is clear that the extensions of the logic $\mathcal{SH}$ form a (complete) lattice.
 Let $\mathbf{L_V(\mathbb{SH})}$ denote the lattice of subvarieties of $\mathbb{SH}$ and $Ext(\mathcal{SH})$ the lattice of axiomatic extensions of the logic $\mathcal{SH}$.
 
The following theorem, due to Blok and Pigozzi, is one of the hallmark accomplishments of Abstract Algebraic Logic (see \cite{BlPi89}).

\begin{Theorem} (\cite[Theorem 3.40]{Fo17}) \label{iso}  Let L be an algebraizable logic with the variety $\mathbb{K}$ as its equivalent algebraic semantics. Then $Ext(\mathcal{L})$ is dually isomorphic to $\mathbf{L_V(\mathbb{SH})}$.
\end{Theorem}

Using Theorem \ref{teorema_SH_implicativa}, Corollary \ref{CorExt} and Theorem \ref{iso},
 the following theorem, crucial to the rest of the paper, is 
 proved in \cite{CoSa22a}.
\begin{Theorem} \label{Ext_subvariety_Theorem}
$Ext(\mathcal{SH})$ is dually isomorphic to $\mathbf{L_V(\mathbb{SH})}$.
\end{Theorem}  

The above theorem justifies the use of the phrase ``the logic corresponding to a subvariety $\mathbb  V$ of $\mathbb{SH}$. 

\smallskip

Let $\mathbf{Mod}(\mathcal E) := \{\mathbf{A} \in \mathbb{SH}: \mathbf{A} \models \delta \approx 1, \mbox{ for every } \delta \in \mathcal E\} $.  

\smallskip
The following theorem is an immediate consequence of Theorem \ref{Ext_subvariety_Theorem} and plays an important role in the rest of the paper.

 \begin{Theorem} \label{completeness_SH_extension} \label{teo_040417_01}
 	Let $\mathcal{E}$ be an axiomatic extension of the logic $\mathcal{SH}$.  Then          
	
     (a)  $\mathcal{E}$ is also algebraizable with the same equivalence formulas and defining equations as those of the logic $\mathcal{SH}$

    (b)  $\mathbf{Mod}(\mathcal E)$ is an equivalent algebraic semantics of $\mathcal E$.   
 \end{Theorem}

 \subsection{The Connexive Semi-Heyting Logic $\mathcal{CSH}$} 

\

\medskip

Let $\neg \alpha := \alpha \to 0$.\\

{\bf Language}:= $\langle \lor, \land, \to, \bot, \top \rangle$\\

{\bf Axioms:}  (A1)-A11) together with the following axioms:\\

\indent  \quad (A12)  \quad $\neg(\neg \alpha \to \alpha)$ \hspace{2.4cm} (AT1) \quad (Aristotle's Thesis 1),
 
\indent \quad (A13)  \quad $\neg(\alpha \to \neg \alpha)$ \hspace{2.4cm} (AT2) \quad (Aristotle's Thesis 2),

\indent \quad (A14)  \quad $(\alpha \to \beta) \to \neg(\alpha \to \neg \beta)$  \quad (BT1) \quad (Boethius' Thesis 1),
 
\indent \quad (A15)  \quad $(\alpha \to \neg \beta) \to \neg(\alpha \to \beta) \ \ \quad (BT2)$ \ \quad (Boethius' Thesis 2).\\

  {\bf Rule of Inference:}\\
 	
\indent \quad (SMP):   From $\phi$ and $\phi \to_H {\gamma}$, deduce $\gamma$  ({\rm semi-Modus Ponens}).\\

It is clear that the logic $\mathcal{CSH}$ is an axiomatic extension of $\mathcal{SH}$. 

We will now introduce some more extensions of the logic $\mathcal{SH}$.

\begin{definition}
\begin{thlist}
\item[a] The logic $\mathcal{AT}i$, $i=1,2$ is an extension of the logic $\mathcal{SH}$ defined by {\rm(ATi)}.  
\item[b] The logic $\mathcal{BT}i$, $i=1,2$ is an extension of the logic $\mathcal{SH}$ defined by {\rm(BTi)}.
\end{thlist}
\end{definition}
Hence we have the following theorem:
 
\begin{Theorem} The logics $\mathcal{CSH}$, $\mathcal{AT}i$, $i=1,2$  and $\mathcal{BT}i$, $i=1,2$ are algebraizable  
with the defining identity $x \approx x \to_H x$ \rm(or equivalently, $x \approx 1$) and the equivalence formulas $\Delta = \{p \to_H q, q \to_H p \}$.

\end{Theorem}

\vspace{.7cm}

\section{Algebraically Speaking: Connexive Semi-Heyting algebras} 

The discussion in the preceding section leads us to the algebraic notion of ``connexive semi-Heyting algebras.''

Consider the following identities: 
      \begin{itemize}
      \item[\rm {i)}] $(x^* \to x)^* \approx 1$ \qquad \qquad \qquad  {\rm(AT1)}  
      \item[{\rm (ii)}] $(x \to  x^*)^* \approx 1$   \qquad \qquad \qquad {\rm(AT2)}      
      \item[{\rm (iii)}] $(x \to y) \to (x \to y^*)^* \approx 1$  \quad {\rm(BT1)}  
      \item[{\rm (iv)}] $(x \to y^*) \to (x \to y)^* \approx 1$  \quad {\rm(BT2)}. 
      \end{itemize}

Note that the names (AT1), etc. just used in the algebraic context may appear to be ambiguous.   However, their use in logic or algebraic context will make their meaning clear.
  
\begin{definition}
A semi-Heyting algebra $\mathbf{A}$ is connexive if $\mathbf{A}$ satisfies the identities \rm(AT1), (AT2), (BT1), and (BT2).
Let $\mathbb{CSH}$ denote the subvariety of $\mathbb{SH}$ consisting of connexive semi-Heyting algebras.
Let $\mathbb{AT}1$, $\mathbb{AT}2$, $\mathbb{BT}1$ and $\mathbb{BT}2$ denote the subvarieties of $\mathbb{SH}$ defined, respectively, by the identities (AT1), (AT2), (BT1), and (BT2).
\end{definition}

The following theorem is immediate from Theorem \ref{completeness_SH_extension}. 
\begin{Theorem}
The varieties $\mathbb{CSH}$, $\mathbb{AT}1$,
$\mathbb{AT}2$, $\mathbb{BT}1$ and $\mathbb{BT}2$ respectively are the equivalent algebraic semantics of the logics $\mathcal{CSH}$, $\mathcal{AT}1$,
$\mathcal{AT}2$, $\mathcal{BT}1$ and $\mathcal{BT}2$.
\end{Theorem}

We note that the algebra $\mathbf{\bar{2}}$ (see Figure 1) and $\mathbf{CSH3}${\footnote
{We warmly thank Professor Paoli for sending us a prepublication copy of \cite{FaLePa}, whose authors have independently observed that $\mathbf{CSH3}$ (named therein as KO3) satisfies connexive logic theses.  The same observation is also made independently in \cite{KaOm21}.}
(which was named $\mathbf{L_9}$ in \cite{Sa07}) (see Figure 2) are connexive, while $\mathbf{L_{10}} \in \mathbb{AT}i$, $i=1,2$.  Note also that the algebra $\mathbf{CSH3} \in \mathbb{CSH}$ fails to satisfy the identity: $(x \to y) \to (y \to x) \approx 1$, since $(a \to 1)\to (1 \to a)=1 \to a=a \neq 1$.
Hence it follows that
the logic $\mathcal{CSH}$ is non-symmetric.

Recently, \cite{FaLePa} has introduced \em{connexive Heyting algebras.  

An algebra $\mathbf A=\langle A, \lor, \land, \to, 0, 1\rangle$, with a bounded distributive lattice-reduct, is a \em{connexive Heyting algebra} if A satisfies:

C1: $(x \to y) \to ((y \to z) \to (x \to z)) \approx 1; $

C2: $(x \to y) \to \neg (x \to \neg y)) \approx 1;$

C3: $x \land (x \to y) \approx x \land y;$

C4: $x \to y \leq (z \land x) \to (z \land y);$

C5: $x \to y \leq (z \lor x)\to (z \lor y);$

Let $\mathbb{CH}$ denote the variety of connexive Heyting algebras. 

It is shown in \cite{FaLePa} that $\mathbb{CH}$ is, in fact, a subvariety of $\mathbb{SH}$. In view of the axiom C2, it follows that $\mathbb{CH}$ is also a subvariety of $\mathbb{CSH}$.
It is also shown in \cite{FaLePa} that $\mathbb{CH}$ is term-equivalent to the variety $\mathbb{H}$ of Heyting algebras.  It, therefore, follows that the latttice of subvarieties of $\mathbb{CH}$ is isomorphic to an interval in that of $\mathbb{BT}1$ and that the latttice of subvarieties of
$\mathbb{BT}1$ has cardinality $2^{\omega}$.

\vspace{1cm}
\section{Mutual Relationships among the varieties: $\mathbb{AT}1$, $\mathbb{AT}2$, $\mathbb{BT}1$ and $\mathbb{BT}2$}

There are some interesting and important connections among these varieties and hence, among their corresponding logics.  In this section, we explore these connections.  

\begin{Theorem} \label{six,one}
Let $\mathbf A \in \mathbb{SH}$ such that $\mathbf A \models  {\rm (AT1)}$. Then $\mathbf A \models {\rm (AT2)}$.   Thus, $\mathbb{AT}1 \subseteq \mathbb{AT}2$.
\end{Theorem}

The proof of the preceding theorem depends on the following lemmas.

Unless stated otherwise, the following hypothesis is assumed in the lemmas below:

``Let $\mathbf A \in \mathbb{SH}$ satisfy the identity {\rm (AT1)}.''

\begin{Lemma} \label{659}   
Let $x \in A$.  Then $x^* \to x  =  0.$  In particular, $0 \to 1=0$.                                                                       
\end{Lemma}

\begin{Proof}
From (AT1),  we have $(x^*  \to  x)^* = 1$, from which it is clear that
 $x^*  \to x  =  0$, since $x^*  \to  x \leq (x^*  \to  x)^{**}$.  Hence, $1^* \to 1=0$, implying $0 \to 1=0$.
\end{Proof}

\begin{Lemma} \label{6113}  
$0 \to x = x^*$.  
\end{Lemma}

\begin{Proof} Let $x \in A$.
 $x \land (0 \to x)= x \land (0 \to 1) = x \land (1^* \to 1) =0$ in view of (SH3) and Lemma \ref{659}.  Thus $ x \land (0 \to x) =0,$
which implies that $0 \to x \leq x^*$.  
Next,
clearly, we have
 $x^* \land (0 \to x) = x^* \land (0 \to 0) = x^*;$   
hence,
 $x^* \leq (0 \to x).$  
Thus we get
 $0 \to x = x^*.$  
\end{Proof}

\bigskip
\begin{Proof}{\bf of Theorem \ref{six,one}}: Let $x \in A$.  Then
$x \land (x \to x^*) =x \land x^*= 0$; whence
$x \to x^* \leq x^*$.  Also, $x^* \land (x \to x^*)=x^* \land (0 \to x^*)  = x^* \land x^{**}=0$ by Lemma \ref{6113}, implying that $x \to x^* \leq x^{**}$.  Thus,
$x \to x^* \leq x^* \land x^{**}=0$, implying $x \to x^*=0$.  It now follows that
 $(x \to x^*)^* =1.$
\end{Proof}

\medskip
The converse holds as well.

\begin{Theorem} \label{AT21} 
Let $\mathbf A \in \mathbb{SH}$ such that $\mathbf A \models  {\rm (AT2)}$.  Then 

$\mathbf A \models {\rm (AT1)}$.   Thus, $\mathbb{AT}2 \subseteq \mathbb{AT}1$.
\end{Theorem}  

The proof of this theorem relies on the following lemma.

\begin{Lemma} \label{LE50}
Let $\mathbf A \in \mathbb{SH}$ satisfy the identity {\rm (AT2)}.  Then
$\mathbf A \models 0 \to x \leq x^*. $  
\end{Lemma}

\begin{Proof} Let $x,y,z \in \mathbf A.$
From (AT2), we have
$(0 \to 0^*)^* = 1,$  
implying\\
 $(0 \to 1)^*  = 1.$  
Hence,
$0 \to 1=0.$  
Therefore, $x \land (0 \to 1)= x \land 0$, which implies
 $x \land ((x \land 0) \to (x \land 1)) = 0,$   
yielding
$x \land (0 \to x) = 0;$ thus, $0 \to x \leq x^*,$ proving the lemma.  
\end{Proof}

\medskip
We are now ready to prove Theorem \ref{AT21}.\\

\begin{Proof}{\bf of Theorem \ref{AT21}:}  Let $x \in A.$  Then
$x \land (x^* \to x)=x \land [(x \land x^*) \to x]=x \land (0 \to x)=0$ by Lemma \ref{LE50}.  Hence,
$x^* \to x \leq x^*$.  Also, $x^* \land (x^* \to x)=x^* \land x=0$, whence $x^* \to x \leq x^{**}.$  Thus,
$x^* \to x \leq x^* \land x^{**}=0$, implying $(x^* \to x)^*=1$, which completes the proof.
\end{Proof}

\begin{Theorem}  
$\mathbb{BT}2 \subset \mathbb{AT}1$.
\end{Theorem}  

\begin{Proof}
Let $\mathbf A \in \mathbb{BT}2$ with $x \in A$. 
From (BT2)
we have  
$(x^* \to x^*) \to (x^* \to x)^* = 1; $  whence,
$1 \to (x^* \to x)^* = 1,$ which 
simplifies to 
$(x^*  \to  x)^* = 1,$ implying $\mathbf A \in \mathbb{AT}1$.  Also,
observe that (AT1) does not imply (BT2), since the semi-Heyting algebra $\mathbf{L_{10}}$ (see Figure 2)
satisfies (AT1) and does not satisfy (BT2) (at 2, 0). 
The proof is complete.
\end{Proof}

\begin{Theorem} \label{BT!2}  $\mathbb{BT}1 \subset \mathbb{BT}2$.
\end{Theorem}

The proof depends on the following lemmas.

Unless stated otherwise, the following hypothesis is assumed in the lemmas below:

``Let $\mathbf A \in \mathbb{SH}$ satisfy the identity {\rm (BT1)}.''

\begin{Lemma} \label{L73} 
Let $x \in A$.  Then
 $(x \to 1) \leq x^{**}.$   
\end{Lemma}

\begin{Proof}
From (BT1), we get
 $(x \to  1) \to  (x \to  1^*)^* = 1,$  
whence,
 $(x \to  1) \to x^{**} = 1.$ 
Hence,
$(x \to  1) \land [(x \to  1) \to x^{**}] =(x \to  1) \land 1,$ implying that
 $(x \to 1) \land x^{**} = x \to 1,$ proving the lemma.
\end{Proof}

\begin{Lemma} \label{L110}  
Let $x, y \in A$.  Then
 $(x \to y^*)^* \leq (x \to y)^{**}$.  
\end{Lemma}

\begin{Proof}
From (BT1), we have
$(x^* \to x) \to (x^* \to x^*)^*=1$, which yields 
 $(x^* \to x) \to 1^* = 1,$ implying (AT1): 
\begin{equation} \label{eqA}
 (x^* \to x)^* = 1.  
 \end{equation}
Again, from (BT1) we have $(x^{**} \to x) \to (x^{**} \to x^*)^*=1$.  Hence, by (\ref{eqA}), we get 
\begin{equation} \label{eqB}
 (x^{**} \to x) \to 1 = 1. 
\end{equation} 
So, 
we have
\begin{align*}
x^{**} &= x^{**} \land 1\\
 &= x^{**} \land [(x^{**} \to x) \to 1] \text{ by (\ref{eqB})}\\  
&=  x^{**} \land [\{x^{**} \land (x^{**} \to x)\} \to (x^{**} \land 1)] \text{ by (SH3)}  \\         
&= x^{**} \land [\{x^{**} \land x)\} \to x^{**}]   \\
&= x^{**} \land (x \to x^{**})  \\
&= x^{**} \land (x \to 1) \text{ by (SH3)}. 
\end{align*}
Thus,
$x^{**} \leq (x \to 1),$  which, together with
 Lemma \ref{L73}, yields
 \begin{equation} \label{L84}
 x^{**} = x \to 1. 
 \end{equation}
Next, Observe
$x \land  (y \to 1)$
 = $x \land  ((x \land  y) \to x) $ 
 = $x \land  (y \to x). $\\   
 Hence, in view of (\ref{L84}), we have
 \begin{equation} \label{L87}
 x \land  y^{**} = x \land  (y \to x).  
 \end{equation}
Replacing $x$ by $(x \to y^*)^*$ and $y$ by $x \to y$ in (\ref{L87}), we get\\
$(x \to y^*)^* \land  (x \to y)^{**} = (x \to y^*)^* \land [(x \to y) \to (x \to y^*)^*],$
Hence by (BT1), we get 
 $(x \to y^*)^* \land  (x \to y)^{**} = (x \to y^*)^* $, which proves the lemma.
\end{Proof}

\begin{Lemma} \label{L111} 
$ (x \to y^*)^* = (x \to y)^{**}. $  
\end{Lemma}

\begin{Proof}
We know
$ (1 \to x) \to (1 \to x^*)^* = 1$ by (BT1); that is, 
$x \to x^{**} = 1.$
From (BT1)
we have 
$(x \to y) \land [(x \to y) \to  (x \to y^*)^* ]= (x \to y) \land 1$, which yields
$(x \to y) \land (x \to y^*)^* = (x \to y) $, whence
$(x \to y)^{**} \land (x \to y^*)^{***} = (x \to y)^{**},$  from which we get
$(x \to y)^{**} \land (x \to y^*)^* = (x \to y)^{**} $; thus, $(x \to y)^{**} \leq (x \to y^*)^*,$
 which, in view of Lemma \ref{L110}, leads to
 $(x \to y^*)^* = (x \to y)^{**}. $ 
\end{Proof}

\medskip
\begin{Proof}{\bf of Theorem \ref{BT!2}:}  \\
Let $\mathbf{A} \in \mathbb{SH}$ such that  $\mathbf{A} \models (BT1)$ and let $x,y \in A$. 
 We know from (BT1) that
 $(x \to y^*) \to (x \to y^{**} )^*= 1. $ 
 Hence by applying Lemma \ref{L111} twice, we have
  $(x \to y^*) \to (x \to y)^{***} = 1. $
whence
 $(x \to y^*) \to (x \to y)^* = 1.$  
Hence $\mathbb{BT}1 \subseteq \mathbb{BT}2$.
To show that the inclusion is strict, use the following $\mathbb{BT}2$-algebra that fails to satisfy (BT1) at (1,2).  It has the universe = $\{0,1,2,3\}$ with $0 <2<3<1$, and $\to$ is as follows:

\medskip
\begin{tabular}{r|rrrr}
$\to$: & 0 & 1 & 2 & 3\\
\hline
    0 & 1 & 0 & 0 & 0 \\
    1 & 0 & 1 & 2 & 3 \\
    2 & 0 & 3 & 1 & 1 \\
    3 & 0 & 1 & 2 & 1
\end{tabular} \hspace{.5cm}\\
The proof is now complete.
\end{Proof} 

\medskip
In view of the preceding results, we have arrived at the main theorem of this section.

\begin{Theorem} \label{T_Rel1} We have:
\begin{itemize}
\item[(1)] $\mathbb{AT}1= \mathbb{AT}2$,

\item[(2)] $\mathbb{BT}1 \subset \mathbb{BT}2  \subset \mathbb{AT}1.$
\end{itemize}
\end{Theorem}

The following corollary is immediate Theorem \ref{T_Rel1}.

\begin{Corollary}
A semi-Heyting algebra is connexive iff it satisfies (BT1).  Thus, $\mathbb{CSH} =\mathbb{BT}1$.\\
\end{Corollary}

\subsection{Logically speaking: Mutual relationships among the theses (AT1), (AT2), (BT1), (BT2)}

\

\medskip
We now express Corollary \ref{T_Rel1} about the relationships among the Aristotle's Theses and Boethius' Theses in the logic $\mathcal{SH}$.

\smallskip
Logics $\mathcal{AT}1$, $\mathcal{AT}2$, $\mathcal{BT}1$ and $\mathcal{BT}2$ are axiomatic extensions of the logic $\mathcal{SH}$ defined by the axioms (AT1), (AT2), (BT1) and (BT2), respectively.

The following corollary follows immediately from Corollary \ref{T_Rel1}.
\begin{Corollary} 
 \begin{thlist}
\item[1] The logics $\mathcal{AT}1$ and $\mathcal{AT}2$ are equivalent.

\item[2] $\mathcal{BT}2 \in Ext(\mathcal{AT}1),$

\item[3] $\mathcal{BT}1 \in Ext(\mathcal{BT}2),$
\end{thlist}
\end{Corollary}

The proof of the above corollary was promised in the paper \cite{CoSa22a}.

\vspace{.7cm}

\section{3-valued Connexive Semi-Heyting algebras}

We mentioned earlier that $\mathbf{CSH3}$ (also denoted by $\mathbf{L_9}$ in \cite{Sa11}) is a connexive semi-Heyting algebra.  Let $\mathbb{CSH}3$ denote the variety of semi-Heyting algebras generated by $\mathbf{CSH3}$.
The purpose of this section is to axiomatize the variety $\mathbb{CSH}3$, relative to $\mathbb{BT}1$ (=$\mathbb{CSH}$). This axiomatization enables us to
(1) present a Hibert-style proof system for the logic $\mathcal{CSH}3$ whose equivalent algebraic semantics is $\mathbb{CSH}3$, and (2) prove that $\mathbb{CSH}3$ is term-equivalent to $\mathbb{V}(\mathbf{H}3)$, the variety generated by the 3-element Heyting algebra $\mathbf{L_1}$ (see Figure 2).
Also, as a by-product, we show that 
the non-symmetric connexive logic $\mathcal{CSH}3$ is deductively equivalent to the the well-known 3-valued intuitionistic logic whose equivalent algebraic semantics is the variety generated by the three element Heyting algebra $\mathbf{L_1}$ (see Figure 2).

We start with some preparatory results. The following theorem, which is needed in the sequel, is immediate from \cite[Theorem 7.5]{Sa07}.

\begin{Theorem}\label{T}
A  subdirectly irreducible semi-Heyting algebra satisfies the following condition:

{\rm(SI)} $x \lor y=1$ implies $x=1$ or $y=1.$
\end{Theorem}

\begin{Lemma}\label{La1}
Let $\mathbf{A} \in \mathbb{BT}1$ and 
let $a \in A$ be a dense element.  Then
$a \to 1 = 1$.
\end{Lemma}

\begin{Proof}
From (BT1) we have 
$(1 \to a)\to (1\to a^*)^*=1$, implying $a \to a^{**}=1$.  Hence $a \to 1=1$, as $a^*=0$ by hypothesis.
\end{Proof}

\begin{Lemma}\label{HeightLemma}
Let $\mathbf{A} \in \mathbb{BT}1$ satisfy:

\rm(a) $x \lor (x \to y) \approx x \lor ((x \to y) \to 1),$   

(b) $x^* \lor (x \to y) \approx (x \lor y) \to y$,  

(c) $x \lor y \approx 1$ implies  $x \approx 1$ or $y \approx 1.$ 

Then the height of $\mathbf{A} \leq 2.$ 
\end{Lemma}

\begin{Proof} Let $x,y, a,b \in A$ such that $0 < a < b$.  We wish to show that $b=1$.  Suppose $b \neq 1$.  Then we wish to arrive at a contradiction. 
From (b), we get
 $x^* \lor (x \to 1) = (x \lor 1) \to 1 =1 \to 1$,  
implying
$ x^* \lor (x \to 1)  = 1. $ 
Hence, by (c), we have
\begin{equation}\label{73}
  x^* = 1 \ \text{  or  } \ x \to 1 =1. 
\end{equation}
Suppose $b \to 1 \neq 1$.  Then  $b^*=1$ by \eqref{73}.  So,
 $a = a \land 1=  a \land b^* = a \land [(a \land b) \to (a \land 0)] =
a \land a^* = 0$, as $a \leq b$;  
whence $a=0$--impossible since $a > 0$.  
Thus, we conclude that 
\begin{equation}\label{79}
b \to 1 = 1.  
\end{equation}
From (BT1) we have
$(b \to 1) \to (b \to 1^*)^*=1.$  Hence, by \eqref{79}, we get
 $1 \to b^{**} = 1$,  
implying
\begin{equation}\label{82}
 b^{**} = 1.  
\end{equation}
From (b), we have
 $b^* \lor (b \to x) = (b \lor x) \to x$; 
hence, in view of \eqref{82}, we get
$b \to x=(b \lor x) \to x$.  So,
 $(b \lor x) \land  (b \to x) = (b \lor x) \land x,$  
implying
\begin{equation}\label{95}
 (b \lor x) \land (b \to x) = x.   
\end{equation}
Next,
suppose $x \to 1 \neq 1$.  Then $ x^*=1$ by \eqref{73}, whence  
 $x = 0$. 
Thus,  $x \to 1 = 1$ or  $x = 0$, which, in turn, implies 
\begin{equation} \label{93}
 x \to 1 = 1 \text{ or } x \land y= 0.
\end{equation}
Replacing $x$ by $b \to x$ and $y$ by $b \lor x$ in \eqref{93}, we have 
$(b \to x) \to 1 = 1$ or $(b \to x) \land (b \lor x) =0$.
So, by \eqref{95}, we get
\begin{equation}\label{EqAB}
(b \to x) \to 1 = 1 \text{ or }  x = 0.   
\end{equation}
Since $a \neq 0$, we get from \eqref{EqAB} that $(b \to a) \to 1=1$.
Hence, by (a), we get $b \lor (b \to a) = b \lor 1$, which simplifies to 
 $b \lor (b \to a) = 1.$  Hence, by (c), $b=1$ or $b \to a=1$, implying $b \to a=1$, as $b \neq 1$ by supposition.  
 So,
 $b \land (b \to a) = b \land 1,$ whence
 $b \land a = b,$ which implies $a=b$,
which contradicts the hypothesis that
 $a < b$. 
This completes the proof. 
\end{Proof}

\begin{Lemma} \label{L0x}
Let $\mathbf{A} \in \mathbb{BT}1$ and $x, y \in A$ such that

\rm(i) $ x^* \lor (x \to y) = (x \lor y) \to y, $  

(ii) $x \lor y = 1$ implies $x = 1$ or $y = 1.$ 
Then

(a) $x \neq 0$ implies $x^* = 0$.

(b) $x \neq 0$ implies $0 \to x =0$.
\end{Lemma}

\begin{Proof}
Let $x \in A$ such that $x>0$.  From (i) we get
$x^* \lor (x \to 1)= (x \lor 1) \to 1$; implying
 $x^* \lor (x \to 1)$  = 1. 
Thus by (ii), we have
 $x^* = 1$ or  $x \to 1 = 1,$   
which implies
 $x=0$ or $x \to 1 = 1.$  
Thus, $x \to 1=1$, as $x \neq 0$ by hypothesis.  Therefore, with $y:=1$, we get from (BT1) that
$1 \to ((x \to 1^*) \to 0) = 1$,  as $\mathbf{A} \in \mathbb{BT}1$,
which implies
 $x^{**}= 1.$  Thus, we have
$x^* = 0,$ proving (a).
Again from (BT1), we have
 $(x^* \to x) \to (x^* \to x^*)^* = 1$, 
simplifying to
$(x^* \to x)^* = 1$. 
Hence
$x^* \to x = 0.$  
So, in view of (a), we have
 $0 \to x = 0,$  
completing the proof.
\end{Proof}

\medskip
\begin{Theorem} 
Let $\mathbb{V}$ be the subvariety of $\mathbb{BT}1$ defined by the following axioms:  

\rm(ii) $x \lor (x \to y)=x \lor ((x \to y) \to 1)$, 

(iii) $x^* \lor (x \to y)=(x \lor y) \to y. $ \\ 
Then $\bar{\mathbf{2}}$ and $\mathbf{CSH3}$ are, up to isomorphism, the only
non-trivial subdirectly irreducible algebras in $\mathbb{V}$. 
\end{Theorem}

\begin{Proof} It is easy to see that $\bar{\mathbf{2}}$ and $\mathbf{CSH3}$ are subdirectly irreducible algebras in $\mathbb{V}$.  For the converse,
let $\mathbf{A}$ be subdirectly irreducible in $\mathbb{V}$ with $|A| >1$.  Then, by Theorem \ref{T}, the condition (c) of Lemma \ref{HeightLemma} holds.  Hence, 
 by Lemma \ref{HeightLemma}, the height of the lattice-reduct of $\mathbf{A}$ is $\leq 2$.  Hence the lattice-reduct of $\mathbf{A}$ is isomorphic to that of $\bar{\mathbf{2}}$ or of $\mathbf{CSH3}$, as every non-zero element in A is dense by Lemma \ref{L0x} (a). 
Regarding the description of the $\to$ operation on A:  First, let us suppose $|A|=2$, and let $A =\{0,1\}$. Then it is obvious that $0 \to 0=1=1 \to 1$ and $1 \to 0=0$.  Also, Lemma \ref{L0x} says that $0 \to 1 =0$, thus $\mathbf{A} \cong \mathbf{\bar{2}}$.  Next, suppose $|A|=3$, and let $A=\{0, a, 1\}$, with $0 < a <1$.  As before, $0 \to 0=1=1 \to 1$, $1 \to 0=0$, and $0 \to 1 =0$.  We also know that $1 \to a=a$.  By Lemma \ref{L0x} we have $0 \to a=0$.
Since $\mathbf{A}$ is pseudocomplemented and $|A|=3$, it follows that $a \to 0=0$.   It is also clear that $a \to a=1$.  Since $a$ is dense,
it follows from Lemma \ref{La1} that $a\to 1=1$.  Hence $\mathbf{A} \cong \mathbf{CSH3}$.  Thus the proof is complete.
\end{Proof}

\medskip
 
 We are now ready to present the main result of this section that gives an axiomatization
 of the variety $\mathbb{CSH}_3$ generated by $\mathbf{CSH3}$ which follows immediately from the previous theorem.

\begin{Corollary} \label{CSH3}
The variety $\mathbb{CSH}_3$ is defined, modulo $\mathbb{BT}1$, by:

\rm(i) $x \lor (x \to y) \approx x \lor ((x \to y) \to 1)$, 

(ii) $x^* \lor (x \to y) \approx (x \lor y) \to y. $  
\end{Corollary}

Interestingly enough, we have the following corollary showing a close relationship between the variety $\mathbb{CSH}_3$ generated by $\mathbf{CSH3}$ and the variety $\mathbb{H}_3$ generated by the 3-element Heyting algebra $\mathbf{L_1}$ (see Figure 2).  \cite{Sa07} has given the following axiomatization for the variety $\mathbb{H}_3$:

\begin{Theorem} \cite[Theorem 11.2]{Sa07} Let $\mathbf{A} \in \mathbb{SH}$.  Then $\mathbf{A} \in \mathbb{H}_3$ iff
$\mathbf{A}$ satisfies the following:

{\rm(I1)} $x^* \lor x^{**}\approx 1$;

{\rm(I2)} $x \lor (x \to y) \approx (x \to y)^* \to x$.

\end{Theorem}

Perhaps, we should mention here that it is now known (but still unpublished) that the condition (I1) is superfluous.

Let $\mathbf{A} =\langle A, \lor, \land, \to, 0, 1\rangle \in \mathbb{CSH}_3$.  Define $h(\mathbf{A}) := \langle A, \lor, \land, \to_H, 0, 1\rangle$, where $\to_H$ is defined by:
$x \to_H y :=x \to (x \land y)$.  Let $\mathbf{E} =\langle E, \lor, \land, \to, 0, 1\rangle \in \mathbb{H}_3$.  Define $c(\mathbf{E}) := \langle E, \lor, \land, \to_C, 0, 1\rangle$, where $\to_C$ is defined by:
$x \to_C y= (x \to y) \land (x^* \to y^*)$.

\begin{Corollary}\label{Cor_termequiv}
The variety $\mathbb{CSH}_3$ is term-equivalent to the variety $\mathbb{H}_3$.  More explicitly,

{\em(i)} \  If $\mathbf{A} \in \mathbb{CSH}_3$, then  $h(\mathbb{A}) \in \mathbb{H}_3$,

{\em(ii)} \ If $\mathbf{E} \in \mathbb{H}_3$, then $c(\mathbb{E}) \in \mathbb{CSH}3$,

 {\em(iii)} For $\mathbf{A} \in \mathbb{CSH}_3$, $c(h(\mathbf A))= \mathbf A$,

 {\em(iv)} For $\mathbf{E} \in \mathbb{H}$, $h(c(\mathbf E))= \mathbf E$.
\end{Corollary}

\begin{Proof}
Since the varieties $\mathbb{CSH}_3$ and $\mathbb{H}_3$ are generated by $\mathbf{CSH}_3$ and $\mathbf{H}_3$ respectively, it suffices to verify all four conditions in these two algebras.  Such verification, being routine, is left to the reader.
\end{Proof}

\subsection{The Connexive Logic $\mathcal{CSH}3$} 

\

\medskip
Let $\alpha \Leftrightarrow_H  \beta$ denote the formula $(\alpha \to_H \beta) \land (\beta \to_H \alpha)$.  

In view of the back-and-forth translation mentioned in Section 3, Corollary \ref{CSH3} allows us to present a Hilbert-style axiomatization for the extension $\mathcal{CSH}3$, relative to the logic $\mathcal{CSH}$ (or equivalently, to $\mathcal{BT}1$). \\

\noindent {\bf AXIOMS OF} $\mathcal{CSH}3$:  Axioms of $\mathcal{BT}1$ PLUS the following axioms:

\rm(i) $(\alpha \lor (\alpha \to \beta)) \Leftrightarrow_H (\alpha \lor ((\alpha \to \beta) \to 1))$, 

(ii) $(\neg\alpha \lor (\alpha \to \beta)) \Leftrightarrow_H ((\alpha \lor \beta) \to \beta). $

\begin{definition}  Let $\mathcal{L}1 = \langle L, \vdash_{\mathcal{L}1} \rangle$ and 
$\mathcal{L}2 = \langle L, \vdash_{\mathcal{L}2} \rangle$ be two logics of language L.   $\mathcal{L}1$ and $\mathcal{L}2$ are {\bf deductively equivalent} if there exist two translations  $\sigma_1$ from $\mathcal{L}1$ to $\mathcal{L}2$  and  $\sigma_2$ from $\mathcal{L}2$ to $\mathcal{L}1$ such that for all $\Gamma \cup \{\phi\}\subseteq Fm_L$,

\rm(i) $\Gamma \vdash_ {L1} \phi$ \ iff  \ $\sigma_1 (\Gamma) \vdash_{\mathcal{L}2} \sigma_1(\phi)$, where $\sigma_1 (\Gamma) :=\{\sigma_1(\psi) : \psi \in \Gamma\}$; 

(ii) $\sigma_1(\sigma_2 (\phi)) \dashv \vdash_{\mathcal{L}2} \phi$.
\end{definition}

\begin{Remark}
In view of Theorem \ref{completeness_SH_extension}, it is clear that the logic  $\mathcal{CSH}3$ is algebraizable with $\mathbb{CSH}3$ as its equivalent algebraic semantics.  As the variety $\mathbb{CSH}3$ is finitely generated,
it is also clear that the logic  
$\mathcal{CSH}3$ is decidable.  It is also not hard to prove that the varieties $\mathbb{V}(\mathbf{\bar{2}})$ and $\mathbb{CSH}3$
have the amalgamation property, using arguments similar to those in \cite{CoSa24}.
\end{Remark}

Let $\mathcal{H}3$ denote the 3-valued intuitionistic logic whose equivalent algebraic semantics is the variety $\mathbb{H}3$ (generated by the 3-element Heyting algebra $\mathbf{H3}$).

\begin{Theorem}
The logics $\mathcal{CSH}3$ and $\mathcal{H}3$ are deductively equivalent to each other.
\end{Theorem}

\begin{Proof}
Consider the translations $\sigma_1$ from $\mathcal{CSH}3$ to $\mathcal{H}3$ and $\sigma_2$ from $\mathcal{H}3$ to $\mathcal{CSH}3$ such that

$\sigma_1(\alpha \to \beta):=  \alpha \to_H \beta$, 

$\sigma_2(\alpha \to \beta):=  (\alpha \to \beta) \land (\alpha^* \to \beta^*),$\\
and the remaining connectives are left unchanged.  Let $\Gamma \cup \{\phi\}\subseteq Fm_L$. Now, suppose
 
$\Gamma \vdash_{\mathcal{CSH}3} \phi$.  Then $\{\psi \approx 1 : \psi \in \Gamma\} \vdash_{\mathbb{CSH}3} \phi \approx 1$, by the implicativeness of  $\mathcal{CSH}3.$
Hence, $\{\sigma_2(\psi) \approx 1 : \psi \in \Gamma\} \vdash_{\mathbb{H}3} \sigma_2(\phi) \approx 1$, by Corollary \ref{Cor_termequiv}.  It follows that
$\sigma_2(\Gamma) \vdash_{\mathcal{H}3} \sigma_2(\phi)$.  The proof of the reverse inequality is similar.
\end{Proof}

\vspace{.7cm}

\section{Commutative Semi-Heyting algebras}

Recall that a semi-Heyting algebra $\mathbf{A}$ is commutative if $\mathbf{A} \models x \to y \approx y \to x$.  The variety of commutative semi-Heyting algebras is denoted by 
$\mathbb{SH}_c$.  Observe that $\bar{\mathbf{2}}$ is both commutative and connexive.  $\bar{\mathbf{2}}$ is an important algebra in our context since 
$\mathbb{V}(\mathbf{\bar{\mathbf{2}}})$ is the only atom in the lattice of subvarieties of $\mathbb{CSH}$--actually in the much larger lattice of subvarieties of the variety $\mathbb{AT}1$.  This is so, because  $\bar{\mathbf{2}}$ is a subalgebra of every member of $\mathbb{AT}1$, since $\mathbb{AT}1 \models 0 \to 1 \approx 0$, as seen in Section 5.

\vspace{.3cm}

\subsection{Anti-Boolean Semi-Heyting algebras} 

\

\medskip
Recall that a semi-Heyting algebra $\mathbf A$ is {\em anti-Boolean} if $\mathbf A \in \mathbb{V}(\mathbf{\bar{\mathbf{2}}})$.  The second author has given several equational bases in \cite[Section 9]{Sa07}.  In this subsection, for the convenience of the reader, we first recall those equational bases and then present a few more new equational bases for $\mathbb{V}(\mathbf{\bar{\mathbf{2}}})$.

\begin{Theorem} \cite[Section 9]{Sa07} \label{antiB}
Each of the following is an equational base for
$\mathbb{V}(\mathbf{\bar{\mathbf{2}}})$, modulo $\mathbb{SH}$:
\begin{thlist}
\item[a] \rm(i) $x \lor x^* \approx 1$,\\
(ii) $0 \to 1 \approx 0$;
             
\item[b]  \rm(i) $x^{**} \approx x$,\\            
 (ii) $0 \to 1 \approx 0$; 
 
 \item[c] \rm(i) $x \to y \leq x^* \lor y, $\\  
 (ii) $0 \to 1 \approx 0$; 
 
  \item[d] \rm(i)  $x \to y \approx x^* \lor y,$\\
 (ii) $0 \to 1 \approx 0$; 
 
  \item[e] \rm(i) $x \lor (y \to z) \approx  (x \lor y) \to (x \lor z),$\\
  (ii) $0 \to 1 \approx 0$; 
  
  \item[f] $x \to (y \to z) \approx (x \to y) \to z$  (associative property of $\to$);
  
  \item[g] $x \approx (x \to y) \to y$.
\end{thlist}
\end{Theorem}

It was also proved in \cite{Sa07} that $\mathbb{V}(\mathbf{\bar{\mathbf{2}}})$ is term-equivalent to the variety of Boolean rings and (hence) to $\mathbb{V}({\mathbf{2}}).$

We will now present new equational bases for $\mathbb{V}(\mathbf{\bar{\mathbf{2}}})$. 

It was shown in \cite{CoSa19} that one of the only three identities of associative type of length 3 that hold in $\mathbb{SH}$ is the exchange property (Ex): $x \to (y \to  z) \approx y \to (x \to z)$ (which was named $\mathcal{A}_4$ there).  The subvariety $\mathbb{EX}$ of $\mathbb{SH}$ defined by (Ex) includes Heyting algebras.
The following theorem shows an interesting interplay between commutativity and (Ex).

\begin{Theorem}
Let $\mathbf A \in \mathbb{SH}$.  Then $\mathbf A$ is an anti-Boolean algebra iff $\mathbf A$ satisfies:

\rm(i) $ x \to y \approx y \to x$ \quad (Commutativity),

(ii) $ x \to (y \to z) \approx  y \to (x \to z) \quad (Ex)$. 
\end{Theorem}

\begin{Proof} Suppose $\mathbf A \in \mathbb{SH}$ satisfies (i) and (ii) and let $x \in A$.  We first prove that $x^{**} = x$.  Now,\\
\begin{align*}
x^{**} &= (x \to 0) \to 0)\\
&= (0 \to (x \to 0)) \text{ by commutativity}\\
&= (x \to (0 \to 0) \text{ by exchange}\\
&=x \to 1 \\
&=1 \to x\\
&=x.
\end{align*}
Thus, $x^{**} = x$. Observe also that $0 \to 1=1 \to 0=0$ by commutativity.
It follows from Theorem \ref{antiB}(b) that $\mathbf A$ is an anti-Boolean algebra.  The converse is trivial.
\end{Proof}

\begin{Remark}
Commutativity does not imply {\rm(Ex)} as shown by the algebra $\mathbf{L_{10}}$ {\rm(see Figure 2)}.  Also, {\rm(Ex)} does not imply commutativity since the variety $\mathbb{EX}$ contains Heyting algebras.  Thus, $\mathbb{EX}$ and $\mathbb{SH}_c$ are incomparable in the lattice of subvarieties of $\mathbb{SH}$.
\end{Remark}

\begin{Theorem}\label{anti-Boolean}
Let $\mathbf A \in \mathbb{SH}$.  Then  $\mathbf A$ is an anti-Boolean algebra iff $\mathbf A$ satisfies:

\rm(i)  $(x^* \to  x)^* \approx 1$  \rm(AT1),

(ii) $x^*  \lor y^* \lor (x \to y) \approx 1. $
\end{Theorem}

\begin{Proof}  Since $\Rightarrow$ is trivial, for the other direction, it suffices, in view of Theorem \ref{antiB}(a), to show that
$\mathbf A$ satisfies $0 \to 1 = 0 $ and $x \lor x^{*} = 1.$
 From (AT1) it immediately follows that $0 \to  1 = 0.$  
Next, with $x:=1$ in (ii), we get
 $0 \lor ((x \to  0) \lor (1 \to x)) = 1$, 
which implies
 $x \lor x^*=1$, completing the proof.
 \end{Proof}

\medskip
 
\begin{Theorem}
Let $\mathbf A \in \mathbb{SH}$.  Then  $\mathbf A$ is an anti-Boolean algebra iff $\mathbf A$ satisfies:

\rm(i)  $(x^* \to  x)^*  \approx 1$  \rm(AT1),

(ii) $ ((x \to  y) \lor (y \to  z)) \lor (z \to  x) \approx 1. $ 
\end{Theorem}

\begin{Proof} Suppose $\mathbf{A}$ satisfies the conditions of the theorem.
We have already seen that (AT1) implies 
 $0 \to  1 = 0.$  
From (ii) we get $(1 \to  x) \lor  ((x \to  0) \lor  (0 \to 1) = 1,$ implying
 $x \lor x^*= 1.$ 
Thus, $\mathbf A$ is an anti-Boolean algebra.  Since the converse is trivially true, the proof is complete.
\end{Proof}

\subsection{Connection between Commutativity and  Connexivity}

 \
\smallskip
It was noted earlier that $\bar{2}$ is a connexive semi-Heyting algebra.  Note that it is also commutative.  Note also that $\mathbf{L_{10}}$ (see Figure 2) is an (AT1)-connexive semi-Heyting algebra that is also commutative.  These remarks lead us naturally to investigate the connection between the commutative identity and the connexive identities.  Recall that $\mathbb{SH}_c$ denotes the variety of commutative semi-Heyting algebras.

In this subsection, we first show that the only commutative semi-Heyting algebras that belong to the variety  $\mathbb{BT}2$ are from $\mathbb{V}(\bar{2})$ and then show that 
$\mathbb{SH}_c \subseteq \mathbb{AT}1.$
 
\begin{Theorem}\label{COMBT2}
$\mathbb{BT}2 \, \cap\ \mathbb{SH}_c = \mathbb{V}(\bar{2})$.  
\end{Theorem}

\begin{Proof}
Let $\mathbf{A} \in \mathbb{BT}2 \cap \mathbb{SH}_c$ and let $x \in A$.
It suffices to prove that $x^{**}=x$ and $0 \to 1 =0$.
From 
 $1 \to x = x$, we get, using commutativity, that
 \begin{equation}\label{E49}
 x = x \to 1.  
 \end{equation}
From (BT2), with $y=0 $, 
we have
 $(x \to (0 \to 0)) \to ((x \to 0) \to 0) = 1,$  
which simplifies to
 $(x \to 1) \to x^{**} = 1,$  
which, by \eqref{E49} and commutativity, gives
 $x^{**} \to x = 1.$  
Hence,
$x^{**} \land (x^{**} \to x) = x^{**} \land 1$, which implies
$x^{**} \land x = x^{**}$.  Hence it follows that $x^{**} = x$.  Also, since $1 \to 0 = 0$, we have by commutativity that $0 \to 1=0.$
This completes the proof.
\end{Proof}

\begin{Corollary}
$\mathbb{BT}1 \, \cap\ \mathbb{SH}_c = \mathbb{V}(\bar{2})$.  
\end{Corollary}

\begin{Proof}
Observe $\bar{2} \in \mathbb{BT}1$  and use Theorem \ref{T_Rel1} and Theorem \ref{COMBT2}.
\end{Proof}

\medskip

In view of Theorem \ref{COMBT2},
the question, as to what  $\mathbb{SH}_c \cap \mathbb{AT}1$ is, naturally arises.
Interestingly, it turns out that commutative semi-Heyting algebras actually form a subvariety of the variety  $\mathbb{AT}1$.  

We need a lemma before we can prove this theorem. 

\begin{Lemma} \label{154}  
 Let $\mathbf A \in \mathbb{SH}_c$ and $x,y \in A$.  Then
  
 $(x \lor y) \land (y \to x) = x \land y$. 
\end{Lemma}

\begin{Proof}
\begin{align*}
(x \lor y) \land (y \to x) \\
&=(x \land (y \to x)) \lor (y \land (y \to x)) \\
&= (x \land (x \to y)) \lor (y \land x) \text{ by commutativity}\\
&=(x \land y) \lor (y \land x)=x \land y,
\end{align*}
proving the lemma.
\end{Proof}

\begin{Theorem}\label{Com-AT1}
Let $\mathbf A \in \mathbb{SH}_c$.  Then $\mathbf A \models \rm(AT1).$ 
Hence, 
$\mathbb{SH}_c \subseteq \mathbb{AT}1.$
\end{Theorem}

\begin{Proof} Let $x,y \in A$.
Observe that
$[(x \lor x^*) \land y]^*=[(x \lor x^*)^{**} \land y^{**}]^*
= [(x^* \land x^{**})^{*} \land y^{**}]^* =(0^* \land y^{**})^* =y^{***}=y^*$
Thus, 
\begin{equation}\label{EqA}
[(x \lor x^*) \land y]^*= y^*.
\end{equation}
Hence, replacing $y$ by $(x^* \to x)^*$, we have
 $[(x \lor x^*) \land (x^* \to x)]^*= (x^* \to x)^*$, which, in view of Lemma \ref{154}, yields $(x \land x^*)^*=(x^* \to x)^*$, whence $ 1= (x^* \to x)^*$.
\end{Proof}

\vspace{.7cm}

\section{Connection between (AT1) and the weak commutative identity: $x^* \to y^* \approx  y^*  \to  x^*$.}

\
Theorem \ref{Com-AT1} motivated us to investigate the identity $x^* \to y^* \approx  y^*  \to  x^*$, which is a slight weakening of the commutative identity.
In this section, we will present an interesting connection between the identities $x^* \to y^* \approx  y^*  \to  x^*$ and (AT1).
Let $\mathbb{S}t\mathbb{SH}$ denote the subvariety of $\mathbb{SH}$ defined by the Stone identity: $x^* \lor x^{**} \approx 1$.
 
\begin{Theorem}\label{ATWEAKCOM} 
Let $\mathbb{S}t\mathbb{SH} \models \rm(AT1)$.   Then 

$\mathbb{S}t\mathbb{SH} \models x^* \to y^* \approx  y^*  \to  x^*$.
\end{Theorem}

The proof of this theorem relies on the following lemma.

\begin{Lemma} \label{LemmaA}
Let $\mathbf A \in \mathbb{SH}$ satisfy \rm(AT1) and let $x \in A$.  Then

$x^*  \to  1 = x^*.$
\end{Lemma}

\begin{Proof}
Since
$x \land ((x \land y) \to z) = x \land (y \to z), $\\  
we have
\begin{equation} \label{58}
 x \land (y \to z) \leq (x \land y) \to z. 
\end{equation}
Hence $x \land (y \to y) \leq (x \land y) \to y$, implying
\begin{equation} \label{70}
 x \leq (x \land y) \to y.  
 \end{equation}
So, substituting $x$ by $x^* \to 1$ and $y$ by $x$ we have
 \begin{align*}
x^* \to 1 &\leq [(x^* \to 1) \land x] \to x\\
                      &= [x \land (0  \to x)] \to x \\
                      &= (x \land x^*) \to x \text{ by Lemma \ref{6113}}\\  
                     &= 0 \to x.  
\end{align*}
Hence by Lemma \ref{6113} we have                           
$x^* \to 1 \leq x^*.$ 
Also, $x^* \land (x^* \to 1)=x^* \land (x^* \to x^*) = x^*$, implying $x^* \leq x^* \to 1$.
Thus we have $x^* \to 1 =x^*,$ 
completing the proof.
\end{Proof}

\bigskip
\begin{Proof}{\bf of Theorem \ref{ATWEAKCOM} }:\\
Let $\mathbf{A} \in \mathbb{S}t\mathbb{SH}$ be subdirectly irreducible with $|A| >1$ and let $x \in A$.
Since $x^* \lor x^{**} = 1$, we have $x^* =1$ or $x^{**}=1,$
implying
 $x^*=1$ or $x^{*}=0$.  So, we consider two cases.  First, suppose $x^*=1$. Then $x^* \to y^* = 1 \to y^* =y^*.$  On the other hand, $y^* \to x^*=y^* \to 1= y^*$ by Lemma \ref{LemmaA}, implying $x^* \to y^* = y^* \to x^*$
Next, suppose $x^*=0$.  Then $x^* \to y^* = 0 \to y^*=y^{**}$ by Lemma \ref{6113};
and $y^* \to x^* =y^* \to 0=y^{**}$. Thus $x^* \to y^*=y^* \to x^*$, proving the theorem.
\end{Proof}

\medskip
We believe that the Stone identity is superfluous in the hypothesis of Theorem \ref{ATWEAKCOM}.  So we state the following conjecture:

\smallskip
{\bf Conjecture}: 
Let $\mathbf A \in \mathbb{SH}$ satisfy \rm(AT1).   Then $\mathbf A \models x^* \to y^* \approx  y^*  \to  x^*$.

\medskip
It turns out that the converse of the above theorem holds as well, even in a more general context.

\begin{Theorem} \label{ConverseAT1} Let $\mathbf A \in \mathbb{SH}$ satisfy :
$x^* \to y^* \approx y^* \to x^*$.  
 Then $\mathbf A \models \rm(AT1).$
\end{Theorem}

The proof of the above theorem depends on the following lemmas.

\begin{Lemma} \label{54} 
Let $\mathbf A \in \mathbb{SH}$ satisfy $x^* \to y^* \approx y^* \to x^*$ and let $x \in A$.  Then 
$x \land (0 \to  x) = 0.$  
\end{Lemma}

\begin{Proof}
From the hypothesis we get
$0^* \to x^*=x^* \to 0^*$, which implies
 $ x^*= x^* \to 1.$ 
Hence, $0 = 0 \to 1$. 
Then we have $x \land 0 = x \land (0 \to 1)$, implying that
 $x \land (0 \to x) = 0.$  by (SH3).
\end{Proof}

\begin{Lemma} \label{958}
Let $\mathbf A \in \mathbb{SH}$ 
and let $x,y \in A$.  Then

 $((x \to y) \to z) \leq [((x \land (y \to z)) \to (y \land ((x \to y) \to z))) \to z].$ 
\end{Lemma}

\begin{Proof}
Replacing x by $(x \to y) \to z)$, y by x, z by y and u by z in Lemma \ref{634}, we get

$((x \to y) \to z) \land (((((x \to y) \to z) \land x) \to (((x \to y) \to z) \land y)) \to z) $
$=(x \to y) \to z) \land ((x \to y) \to z)= (x \to y) \to z)$
Hence, by   
Lemma \ref{634}(6)
we have
 $((x \to y) \to z) \land (((x \land (y \to z)) \to (((x \to y) \to z) \land y)) \to z) = (x \to y) \to z. $ 
\end{Proof}

\begin{Lemma} \label{Lem764}
Let $\mathbf A \in \mathbb{SH}$ 
and let $x,y \in A$.  Then

 $((x \to  y) \to  z) \land  (((x \land  (y \to  z)) \to  y) \to  z) = (x \to  y) \to  z.$  
\end{Lemma}

\begin{Proof}
\begin{align*}
x \land ((y \to z) \to u) \\
&= x \land (((x \land y) \to (x \land z)) \to u) &\text{ by Lemma \ref{634}} \\ 
&= x \land  (((x \land  y) \to  z) \to  u) &\text{ by Lemma \ref{634}} \ 
\end{align*} 
Therefore, 
we obtain
\begin{equation}\label{x61}
  x \land  ((y \to  z) \to  u) = x \land  (((x \land  y) \to  z) \to  u). 
\end{equation}
Hence by Lemma \ref{634}  
and \eqref{x61} we have
 \[x \land  ((y \to  (z \land  x)) \to  u) = x \land  ((y \to  z) \to  u).\]  
Hence by Lemma \ref{958}              
we have
 \[((x \to  y) \to  z) \land  (((x \land  (y \to  z)) \to  y) \to  z) = (x \to  y) \to  z,\]  
proving the lemma.
\end{Proof}

\bigskip
\begin{Proof}{\bf of Theorem \ref{ConverseAT1}:}\\
Replacing $y$ by $0$ and $z$ by $x$ in Lemma \ref{Lem764} we get\\
$((x^* \to  x) \land  (((x \land (0 \to x)) \to  0) \to  x) =x^* \to  x $. 
Hence, we have
\begin{align*}
 x^* \to  x 
&=((x^* \to  x) \land  ((0 \to  0) \to  x)  & \text{by Lemma \ref{54}} \\      
 &=(x^* \to  x) \land x \\ 
&= x \land  (x^* \to  x)  \\   
 &=x \land  (0 \to  x) &\text{ by (SH2)} \\ 
 &=0 &\text{by Lemma \ref{54}}.   
 \end{align*}
 Thus $x^* \to  x = 0,$  which implies $(x^* \to  x)^* = 1.$
This completes the proof of the theorem.
\end{Proof}

\begin{Corollary}
Let $\mathbf{A} \in \mathbb{S}t\mathbb{SH}$. Then the following are equivalent:

[\rm(i) $\mathbf{A} \models {\rm(AT1)}$,

(ii) $\mathbf{A} \models x^* \to y^* \approx y^* \to x^*$.
\end{Corollary}

\begin{Proof}
Use Theorem \ref{ATWEAKCOM} and Theorem \ref{ConverseAT1}.
\end{Proof}

\vspace{.7cm}

\section{Connection between (AT1) and  $0 \to 1 \approx 0$ }
\smallskip
Recall that a semi-Heyting algebra $\mathbf A$ is an {\em anti-Heyting algebra} if $\mathbf A \models 0 \to 1 \approx 0$.  The variety of anti-Heyting algebras is denoted by $\mathbb{AH}$. Note that anti-Heyting algebras generalize anti-Boolean algebras.
The purpose of this section is to show that $\mathbb{AH} = \mathbb{AT}1$.  

\begin{Lemma}\label{61}
Let $\mathbf{A} \in \mathbb{AH}$.  Then $\mathbf{A}$ satisfies
 $(0 \to x) \leq x^*.$   
\end{Lemma}

\begin{Proof} 
Since $0 \to 1=0$ we have  
$ x \land (0 \to x) =x \land (0 \to1)=x \land 0=0$, proving the lemma.
\end{Proof}

\begin{Lemma} \label{172}
Let $\mathbf{A} \in \mathbb{AH}$ Then
 $\mathbf{A} \models 0 \to x \approx x^*.$   
\end{Lemma}

\begin{Proof} 
 Substituting $x$ by $x^*$, $y$ by $0$ and $z$ by $x$ in the identity 
$x \land [y \to (x \land z)] \approx x \land (y \to z)$ (which holds in $\mathbf{A}$), 
we obtain
$x^* \land [0 \to (x^* \land x)] = x^* \land (0 \to x),$ which, by Lemma \ref{61},
leads to 
 $x^* \land (0 \to 0) = 0 \to x.$ 
 Hence,
 $x^* = 0 \to x.$  
 \end{Proof}

\smallskip
Now we are ready to prove our main result of this section.

\begin{Theorem} \label{ASHAT}
Let A be a semi-Heyting algebra.  Then\\
$0 \to 1 \approx 0$ is equivalent to {\rm(AT1)}.  Thus $\mathbb{AH} = \mathbb{AT}1$.
\end{Theorem}

\begin{Proof}  
 See Lemma \ref{659}  for the proof that (AT1) implies $0 \to 1 = 0.$
 For the converse let us observe that, by Lemma \ref{172}, we have that $0 \to x = x^*$. Hence, using (SH3),  
 \begin{align*}
 x^* \to x &= (x^* \to x) \wedge (x^* \to x) \\
 &= (x^* \to x) \wedge ((x^*\wedge x) \to x) \\
 &= (x^* \to x) \wedge (0 \to x) \\
 &= (x^* \to x) \wedge x^* \\
 &= x^*\wedge x = 0. 
 \end{align*}
 Hence $(x^* \to x)^* =1$, completing the proof.   
\end{Proof}

\vspace{.7cm}

\section{The varieties $\mathbb{AT}1$, $\mathbb{EX}$ and  $\mathbb{BT}1$}

It is interesting to see that (Ex) has an interesting connection with (BT1).  In this section, we show that (AT1) and (Ex) imply (BT1). 

\begin{Theorem} \label{ThBT1}
Let $\mathbf{A} \in \mathbb{SH}$ satisfy: 

\rm(i) $(x^* \to x)^* \approx 1$  (AT1),

(ii) $x \to (y \to  z) \approx y \to (x \to z)$  (Ex).

\noindent Then (BT1) holds in $\mathbf{A}$.  Thus, $\mathbb{AT}1 \cap \mathbb{EX} \ \subset \mathbb{BT}1.$
\end{Theorem}

We need some preliminary results before we can prove this theorem. 

\begin{Lemma}\label{33}
Let $\mathbf A$ satisfy \rm(AT1) and (Ex) and let $x \in A$.  Then

 $x^* = 0 \to (x \to 1).$ 
\end{Lemma}

\begin{Proof}
We have by Lemma \ref{659} that 
$ x^* \to x = 0.$  
In particular, we have
 $0 \to 1 = 0.$  
From (Ex) we get $x \to 0 = x \to (0 \to 1) = 0 \to (x \to 1)$; hence
 $x^* = 0 \to (x \to 1).$  
\end{Proof}

\begin{Lemma}\label{xx1:12}
Let $\mathbf A$ satisfy \rm(AT1) and (Ex) and let $x, y \in A$.  Then
 $x \to  ((x \to  y) \to  y) = 1.$
\end{Lemma}

\begin{Proof}
$x \to  ((x \to  y) \to  y) =(x \to y) \to (x \to y)=1$
\end{Proof}

\begin{Lemma}\label{41}
Let $\mathbf A$ satisfy {\rm(AT1)} and (Ex), and let $x,y \in A$  Then

\rm(a) $y \to (x \to y)=x \to 1$,  

(b)  $x \to (0 \to x)=0$,  

(c) $(0 \to x)^* = x \to 1.$
\end{Lemma}

\begin{Proof}
From (Ex), we have $y \to  (x \to  y) = x \to (y \to  y) = x \to 1,$  proving (a).
 We have by Lemma \ref{659} 
$0 \to 1 = 0.$  Hence, we have by (Ex) that $x \to (0 \to x) =0 \to (x \to x)
= 0 \to 1=0$; thus, $x \to (0 \to x)= 0,$ proving (b).
 Substituting $y$ by $0 \to x$ in (a) we get $(0 \to x) \to [x \to (0 \to x)]=x \to 1$, which, in view of (b), implies
  $(0 \to x) \to 0 = x \to 1,$ proving (c).  
\end{Proof}

\begin{Lemma}\label{50} 
Let $\mathbf A$ satisfy \rm(AT1) and (Ex) and let $x \in A$.  Then
$x^* =  (x \to 1)^*.$  
\end{Lemma}

\begin{Proof}
Since $x^* \to x =0$, we get from Lemma \ref{41} (a) that
 $x \to (x^* \to x) = x^* \to 1,$ whence
 $x^* = x \to (x^* \to x) = x^* \to 1$. 
Thus 
\begin{equation} \label{35}
x^* = x^* \to 1.   
\end{equation}
Also, $x \to 1 = x^* \to (x \to x^*)$ by Lemma \ref{41} (a), which, as $x \to x^*=0$, implies
\begin{equation} \label{1000}
   x \to 1 =x^{**}.  
 \end{equation} 
Hence we have
\begin{align*}
 (x \to 1)^*  &=x^{***}  &\text{ by \eqref{1000}}\\
 &= x^* \\
 &=x^* \to 1 &\text{ by \eqref{35}},
\end{align*}
proving the lemma.
\end{Proof}

\begin{Lemma}\label{xx1:54} 
Let $\mathbf A$ satisfy \rm(AT1) and (Ex) and let $x \in A$.  Then
$x \to (x \to 1) = 1$.
\end{Lemma}

\begin{Proof}
Substituting $x$ by $0 \to x$ and $y$ by $0$ in Lemma \ref{xx1:12}, we obtain that
$(0 \to x) \to [((0 \to x) \to 0) \to 0]=1$, which, by Lemma \ref{41}(a), implies
 $(0 \to  x) \to  [(x \to  1) \to  0] = 1,$  
which, by Lemma \ref{50}, yields
 $(0 \to  x) \to  (x \to  0) = 1. $  
Hence, $x \to  ((0 \to  x) \to  0) = 1$ by (Ex),  
which, in view of Lemma \ref{41}(c), leads to
 $x \to  (x \to  1) = 1.$  
\end{Proof}

\begin{Lemma}\label{crucial} 
Let $\mathbf{A}$ satisfy (AT1) and (Ex) with $a, b \in A$.  Then\\
$ a \to (b \to 1) = (a  \to  b) \to 1.$
\end{Lemma}

\begin{Proof}
From 
 Lemma \ref{xx1:54} we have $y \to (z \to (y \to z))=1$. Hence by Lemma \ref{41}(a) we get
 $x \to 1=x \to (y \to (z \to (y \to z)))$,
which, by (Ex), yields
 $x \to 1 = y \to (x \to (z \to (y \to z)))$,   
from which, after replacing $x$ by $a \to b$, $y$ by $a$ and $z$ by $b$, we get
$(a \to b) \to 1 = a \to [(a \to b) \to (b \to (a \to b))))]
= a \to (b \to 1)$ by Lemma \ref{41} (a).
\end{Proof}

\medskip
\begin{Proof}{\bf of Theorem \ref{ThBT1}}:
Let $a,b \in A.$  Then
\begin{align*}
 (a \to b) \to (a \to b^*)^* \\                     
 &\hspace{-2.5cm} =(a \to b) \to [\{a \to (b \to 0)\} \to 0]\\   
 &\hspace{-2.5cm}=(a \to  b) \to  [\{a \to  (0 \to  (b \to  1))\} \to  0] \text{ by Lemma \ref{33}}\\ 
&\hspace{-2.5cm}= (a \to  b) \to  [\{0 \to  (a \to  (b \to  1))\} \to  0]  \text{ by (Ex)} \\ &\hspace{-2.5cm}=(a \to b) \to [\{a \to (b \to 1)\} \to 1]  \text{ by Lemma \ref{41}(c)} \\ 
&\hspace{-2.5cm}=(a \to b) \to [\{(a \to b) \to 1\} \to 1]  \text{ by Lemma \ref{crucial}}\\ 
&\hspace{-2.5cm}=\{(a \to b) \to 1\} \to [(a \to b) \to 1]  \text{ by (Ex)}\\
&\hspace{-2.5cm}= 1,
\end{align*}
proving the theorem.
\end{Proof}

\begin{Remark}
The converse of the preceding theorem fails since
{\rm(BT1)} does not imply {\rm (Ex)}: $x \to (y \to z) \approx y \to (x \to z)$. For, the following semi-Heyting algebra, with the chain $0<2<4<3<1$ as its lattice-reduct, satisfies {\rm(BT1)} but  does not satisfy {\rm (Ex)} {\rm(at $2,3,4$)}:\\

\begin{tabular}{r|rrrrr}
$\to$: & $0$ & $1$ & $2$ & $3$ & $4$\\
\hline
    $0$ & $1$ & $0$ & $0$ & $0$ & $0$ \\
    $1$ & $0$ & $1$ & $2$ & $3$ & $4$ \\
    $2$ & $0$ & $1$ & $1$ & $1$ & $3$ \\
    $3$ & $0$ & $1$ & $2$ & $1$ & $4$ \\
    $4$ & $0$ & $1$ & $2$ & $1$ & $1$
\end{tabular} \hspace{.5cm}
\begin{tabular}{r|rrrrr}
$\land$: & $0$ & $1$ & $2$ & $3$ & $4$\\
\hline
   $ 0$ & $0$ & $0$ & $0$ & $0$ & $0$ \\
    $1$ & $0$ & $1$ & $2$ & $3$ & $4$ \\
    $2$ & $0$ & $2$ & $2$ & $2$ & $2$ \\
    $3$ & $0$ & $3$ & $2$ & $3$ & $4$ \\
    $4$ & $0$ & $4$ & $2$ & $4$ & $4$
\end{tabular} \hspace{.5cm}
\begin{tabular}{r|rrrrr}
$\lor$: & 0 & 1 & 2 & 3 & 4\\
\hline
    $0 $& $0$ & $1$ & $2$ & $3$ & $4$ \\
    $1$ & $1$ & $1$ & $1$ & $1$ & $1$ \\
    $2$ & $2$ & $1$ & $2$ & $3$ & $4$ \\
    $3$ & $3$ & $1$ & $3$ & $3 $ & $3$ \\
    $4$ & $4$ & $1$ & $4$ & $3$ & $4$
\end{tabular}

\medskip
This example also shows that {\rm(BT2)} or {\rm(AT1)} does not imply {\rm(Ex)}.
Also, the Boolean algebra $\mathbf{2}$ shows that
{\rm(Ex)} does not imply {\rm(BT1)}.   Thus $\mathbb{EX}$ is incomparable with $\mathbb{BT}1$, $\mathbb{BT}2$ and $\mathbb{AT}1$ in the lattice of subvarieties of $\mathbb{SH}$.
\end{Remark}

\medskip

\section{Concluding Remarks}

In this paper, we have initiated the investigation of the structure  of the lattice of subvarieties of $\mathbb{BT}1$ (=$\mathbb{CSH}$). 

The results in Section 5 of this paper clearly show that $0 \to 1 \approx 0$ holds in all members of $\mathbb{BT}1$.  So, the algebra $\bar{\mathbf{2}}$ is a subalgebra of every member of $\mathbb{BT}1$.
It follows that the variety $\mathbb{AB}$ of anti-Boolean algebras is the only atom in $\mathbf{L}_{\mathbf{V}} (\mathbb{BT}1)$.   Actually, $\bar{\mathbf{2}}$ is a subalgebra of every member of $\mathbb{AT}1$.  So, $\mathbb{AB}$ is the only atom even in $\mathbf{L}_{\mathbf{V}} (\mathbb{AT}1)$.
A similar argument shows that this atom has exactly one cover, namely the variety $\mathbb{CSH}3$, in $\mathbf{L}_{\mathbf{V}} (\mathbb{BT}1)$.   However, the atom $\mathbb{AB}$ of  $\mathbf{L}_{\mathbf{V}} (\mathbb{AT}1)$ has two covers: $\mathbb{CSH}3$ and $\mathbb{V}(\mathbf{L_{10}})$ generated by the 3-element commutative semi-Heyting algebra $\mathbf{L_{10}}$.  

These observations can be easily translated into Ext($\mathcal{CSH}$) using the back-and-forth procedure given in Section 3.  It is also not hard to see that all logics in Ext($\mathcal{CSH}$), except the logic $\mathcal{AB}$, are strongly connexive.  Thus, there is an ample supply of strongly connexive logics--$2^{\omega}$-many of them! 

The preceding observations naturally lead us to raise the following problem for future research:\\

{\bf PROBLEM 1:} Investigate the structure of the lattice of subvarieties of $\mathbb{BT}1$; similarly,
for $\mathbb{BT}2$ and $\mathbb{AT}1$.\\

{\bf PROBLEM 2:} Investigate the ``intermediate logics'' between $\mathcal{CSH}$ and $\mathcal{AB}.$

Since amalgamation property has connections with logical property, the following problem is of interest:\\

{\bf PROBLEM 3:} Describe the subvarieties of the variety $\mathbb{CSH}$ that have the amalgamation property.\\

More generally, the following project would be of great interest.\\

{\bf PROJECT:} Investigate the subvarieties of the variety $\mathbb{CSH}$ with respect to the known variations of amalgamation property, such as strong amalgamation property, super amalgamation property, etc.\\

 The following conjecture has already mentioned in Section 8. \\

{\bf Conjecture}: 
Let $\mathbf A \in \mathbb{SH}$ satisfy \rm(AT1).   Then $\mathbf A \models x^* \to y^* \approx  y^*  \to  x^*$.\\

{\bf PROPOSALS:}

Recall from Section 5 that each of the identity ``$0 \to 1 \approx 0$''  is equivalent to the identity (AT1), and  (BT1) implies (AT1).  Therefore, they 
hold in all of $2^{\aleph_0}$ varieties of non-symmetric, connexive semi-Heyting algebra (equivalently, of $\mathbb{BT}1$-algebras).  Furthermore, (AT1) is already considered a thesis of Connexive Logic.  Based on all these reasons, it is but natural that we make, philosophically speaking, the following proposals: \\ 

{\bf PROPOSAL 1}: The formula ``$0 \to 1 \leftrightarrow 0$'' is a thesis of Connexive Logic.

{\bf PROPOSAL 2}: The formula ``$\neg \phi \to \neg \psi \leftrightarrow \neg \psi \to \neg \phi$'' is a thesis of Connexive Logic.

\end{document}